\documentclass[11pt]{article}
\usepackage{amsmath,amsfonts,amsthm,amssymb, graphicx}

\theoremstyle{plain}
\newtheorem{thm}{Theorem}
\newtheorem{proposition}{Proposition}
\newtheorem{lemma}{Lemma}

\newtheorem{assumption}{Assumption}

\input{tcilatex}
\begin{document}

\title{Corrections to the Central Limit Theorem for Heavy-tailed Probability Densities}
\author{Henry Lam and Jose Blanchet \\
{\it Harvard University and Columbia University} \\
\ \\
Damian Burch and Martin Bazant\\
{\it Massachusetts Institute of Technology} }
\date{} \maketitle

\begin{abstract}
Classical Edgeworth expansions provide asymptotic correction terms
to the Central Limit Theorem (CLT) up to an order that depends on
the number of moments available. In this paper, we provide
subsequent correction terms beyond those given by a standard
Edgeworth expansion in the general case of regularly varying
distributions with diverging moments (beyond the second). The
subsequent terms can be expressed in a simple closed form in terms
of certain special functions (Dawson's integral and parabolic
cylinder functions), and there are qualitative differences depending
on whether the number of moments available is even, odd or not an
integer, and whether the distributions are symmetric or not. If the
increments have an even number of moments, then additional
logarithmic corrections must also be incorporated in the expansion
parameter. An interesting feature of our correction terms for the
CLT is that they become dominant outside the central region and
blend naturally with known large-deviation asymptotics when these
are applied formally to the spatial scales of the CLT.
\end{abstract}

\section{Introduction}

Let $\left( X_{k}:k\geq 1\right) $ be a sequence of independent and
identically distributed (iid) random variables (rv's) with mean zero
and finite variance $\sigma ^{2}$. If $S_{n}=X_{1}+...+X_{n}$, then
the standard
Central Limit Theorem (CLT) says that for each fixed $x\in \mathbb{R}$,%
\begin{equation*}
P\left( S_{n}/n^{1/2}\leq x\right) =P\left( N\left( 0,\sigma
^{2}\right) \leq x\right) +o\left( 1\right)
\end{equation*}%
as $n\nearrow \infty $. Edgeworth expansions provide
additional information about the previous error term, but only if higher moments exist. In particular, if
$E\left\vert X_{k}\right\vert ^{m}<\infty $ for $m\geq 3$ and
assuming that the $X_{k}$'s are continuous random variables (or
simply non-lattice) an Edgeworth
expansion of order $\left\lfloor m-2\right\rfloor $ takes the form%
\begin{equation}
P\left( S_{n}/n^{1/2}\leq x\right) =P\left( N\left( 0,\sigma
^{2}\right) \leq x\right) +\sum_{j=1}^{\left\lfloor m-2\right\rfloor
}n^{-j/2}P_{j}\left( x\right) +o\left( n^{-\left\lfloor
m-2\right\rfloor /2}\right)  \label{Ed1}
\end{equation}%
as $n\nearrow \infty $, where $P_{j}\left( \cdot \right) $ is
appropriately defined in terms of Hermite polynomials and involves
coefficients that depend on $EX_{k}^{j+2}$, for fixed $x\in\mathbb{R}$ (see Section 2 below). As
a consequence, a traditional Edgeworth expansion cannot continue
beyond the $\left\lfloor m-2\right\rfloor $-th correction term in a
situation where $E\left( \left\vert X_{k}\right\vert ^{\left\lfloor
m-2\right\rfloor +3}\right) =\infty $. Our focus here is precisely
on this situation, building on the preliminary results of Bazant (2006).

In particular, our goal in this paper is to
provide additional correction terms to the
Edgeworth expansion and further reduce the order of the error term as $%
n\nearrow \infty $ in the presence of regularly varying
distributions (which are basically power laws; see Equation
(\ref{density})). As we shall discuss in detail below, an
interesting feature of our development is the fact that the
additional correction terms are qualitatively different (giving rise
to expressions involving $\log \left( n\right) $ in some cases)
depending on whether the power-law decaying rate of the distribution
is an integer number or not, or, when such rate is an integer, on
its parity.

After this paper was completed we learned that a related set of
expansion is given in Chapter 2 of Vinogradov (1994). Our results
here appear to be more general than those expansions. In particular,
our results apply to regularly varying tail of any real number $\alpha>2$, while
Vinogradov considers only power-law tail with integer $\alpha$.
Also, our analysis allows us to obtain explicit expressions for the
coefficients in our Edgeworth expansion in terms of special
functions, an important interest in this paper that was not pursued
by Vinogradov.

Our results can be applied in the context of random walks (Hughes
(1995); Weiss (1994)), where corrections to the CLT are sometimes called
Gram-Charlier expansions (and extended to vector random variables).
As is typical for such expansions, our asymptotic formulae can be
very accurate for small $n$, but unlike standard Gram-Charlier
expansions, our corrections also become dominant at the edge of the
central region and beyond, as illustrated by numerical examples in
Bazant (2006) (see also Theorem \ref{moderatedeviations} in Section
4). Random walks with heavy-tailed distributions (or synonymously
``fat-tailed distributions" in physics and other literatures), in
particular those that are regularly varying, are increasingly used
to describe stochastic processes in diverse fields, ranging from
physics (Bouchaud and Georges (1990); Metzler and Klafter (2000)) to
economics (Bouchaud and Potters (2000); Embrechts, Kluppelberg and
Mikosch (1997)).

Edgeworth expansions are also often used in statistics and applied
probability. For instance, in statistics these types of expansions
have enjoyed successful applications in statistical analysis
involving relatively small samples (Field and Ronchetti (1990)). In
bootstrapping, where data are resampled to calculate quantities
related to the estimator, Edgeworth expansion can reduce the order
of the approximation error (Hall (1992)). In applied probability
settings such as queueing or insurance one is often interested in
approximating performance measures of a stochastic system, for
instance the tail of the delay in a queueing system or the ruin
probability in insurance, by means of a CLT. In some cases, one can
derive corrections to such approximations (Blanchet and Glynn (2006)). These corrections are often called corrected diffusion
approximations and are similar in spirit to Edgeworth expansions. In
fact, our ideas here can be applied to developing higher order
corrections to a class of queueing models in heavy-traffic; such
development will be explored in future work.

Other types of approximation that are often applied to the distribution of $%
S_{n}$ include so-called large deviations or moderate deviations
asymptotics. In contrast to Edgeworth expansions, which approximate
the distribution of $S_{n}$ in spatial scales of order $O\left(
n^{1/2}\right) $, large deviations results typically provide
approximations in spatial scales of order $O\left( n\right) $,
whereas moderate deviations are derived in
spatial scales of order $O\left( n^{1/2+\varepsilon }\right) $ for $%
\varepsilon \in (0,1/2)$. In situations where the increment
distribution has finite logarithmic moment generating function,
standard results provide an exponential decay rate to the tail of
$S_{n}$ in large deviations or moderate deviations scales, say for
quantiles that are of the form $nx$. Using the fact that typically
the rate function is twice continuously differentiable at the origin
(as a function of $x$), one can relate (or transition from) large,
moderate and CLT approximations smoothly. On the other hand, when
the increment distributions are regularly varying, the large
deviations rate, as we shall recall in future sections, turns out to
be basically polynomial. Given that further terms in an Edgeworth
expansion reduce the polynomial (in powers of $n^{1/2}$) error rate
of the tail estimate, one may wonder how Edgeworth expansions relate
to large deviations rates as additional correction terms are added.
As we shall see, perhaps surprisingly, our additional terms blend
smoothly with large and moderate deviations rates when these rates
are formally applied to spatial scales of order $n^{1/2}$ only when
the power of the regularly varying tail is non-even. These formal
connections to large deviations asymptotics are further explored in
Section 4; there we would discuss our asymptotics in relation to the
classical uniform large deviations result by Rozovskii (1989), in
which the sum of iid regularly varying rv can be approximated by a
sum of Gaussian distribution and heavy-tail asymptotic. We will also
show that under restrictive conditions our asymptotic provides a
refinement of this result in the moderate deviations region.

As we indicated earlier, our development gives rise to qualitative
differences in the form of the correction terms depending on the
decay rate at infinity of the increment distributions. In order to
further explain such qualitative differences, let us introduce some
notation. We impose the assumption that the $%
X_{k}$'s have a positive density $f_{X}\left( \cdot \right) $ of the form%
\begin{equation*}
f_{X}\left( x\right) =\left( 1+\left\vert x\right\vert \right)
^{-\left( \alpha +1\right) }L\left( \left\vert x\right\vert \right)
, \label{density}
\end{equation*}%
where $L\left( \cdot \right) $ is a slowly varying function (i.e.
$L\left( xt\right) /L\left( x\right) \longrightarrow 1$ as
$x\nearrow \infty $ for every $t>0$) so that $var(X_k)=1$ (our development also deals
with the non-symmetric case; see Section 2.2 and 5). Also assume that $\int_{-\infty}^\infty|\phi(\theta)|d\theta<\infty$ where $\phi(\theta)=Ee^{i\theta X_k}$ is the characteristic function. Denote $\eta(x)$ as the standard
normal density. Theorem \ref{symmetrictheorem} in Section 2 yields
that the density of $S_n/n^{1/2}$ behaves asymptotically as
\begin{equation*}
f_{S_n/n^{1/2}}(x)= \eta(x)\left[1+\sum_{\substack{ 3\leq j<\alpha  \\ j\text{\ even}}}%
\frac{G_j(x)}{n^{j/2-1}}\right]+\left\{\begin{array}{ll}\frac{L(n^{1/2}x)}{n^{\alpha/2-1}}G_\alpha(x)+o(\frac{L(n^{1/2}x)}{n^{\alpha/2-1}})&\text{when $\alpha$ is not even}\\
\frac{\zeta_L(n^{1/2}x)}{n^{\alpha/2-1}}G_\alpha(x)+o(\frac{\zeta_L(n^{1/2}x)}{n^{\alpha/2-1}})&\text{when
$\alpha$ is even}\end{array}\right.
\end{equation*}
as $n\nearrow\infty$, where $\zeta_L(x)$ is a suitable slowly
varying function (see Equation (\ref{zetadefinition})) and
$G_\alpha(x)$ is written in terms of special functions whose
qualitative behaviors depend on the parity of $\alpha$ (see Equation
(\ref{coefficient})).

Our approach to developing the Edgeworth expansions is standard.
First we develop an appropriate Gram-Charlier series for the Fourier
transform and then apply inversion. An interesting element of the
analysis in our situation, however, involves the contribution of a
non-analytic component in the Fourier transform. We shall build the
elements of our expansion starting from the symmetric case and then
extend the construction of our approximation to more general
regularly varying densities. Finally, we note that our approach can
be used to develop asymptotic expansion for stable limits, which
will be reported in future research.

The paper is organized as follows. First, we discuss our main
results in Section 2. Then we discuss in detail the symmetric case
in Section 3. The formal connection to large deviations is given in
Section 4. Section 5 studies the case of non-symmetric density.

\bigskip

\section{Main Results}

\subsection{Review of Ordinary Edgeworth Expansion}

Let $\left( X_{k}:k\geq 1\right) $ be a sequence of iid rv's with
mean zero, unit variance and with a regularly varying density. More
precisely, we
assume that the density $f\left( \cdot \right) $ of $X_{k}$ satisfies%
\begin{eqnarray}
f\left( x\right) &=&\left( 1+x\right) ^{-\left( 1+\beta \right)
}L_{+}\left(
x\right) I\left( x\geq 0\right) \notag\\
&&+\left( 1-x\right) ^{-\left( 1+\gamma \right) }L_{-}\left(
-x\right) I\left( x<0\right) , \label{density}
\end{eqnarray}%
where $\beta ,\gamma >2$ and $L_{+}\left( \cdot \right) $,
$L_{-}\left( \cdot \right) $ are slowly varying functions at
infinity (i.e. $L_{\pm
}\left( cb\right) /L_{\pm }\left( b\right) \longrightarrow 1$ as $%
b\longrightarrow \infty $) defined on $\mathbb{R}$. We shall write $\phi \left( \theta
\right)
=E\exp \left( i\theta X_{k}\right) $ for the characteristic function of $%
X_{k}$ and put $\psi \left( \theta \right) =\log \phi \left( \theta \right) $%
. Assume $\int_{-\infty}^\infty|\phi(\theta)|d\theta<\infty$. Finally, we let $S_{n}=X_{1}+...+X_{n}$, $n\geq 1$, be the random
walk generated by the $X_{k}$'s, set $\phi _{n}\left( \theta \right)
=E\exp \left( \theta S_{n}/n^{1/2}\right) =\phi^{n}\left( \theta
/n^{1/2}\right)$ and define $\psi _{n}\left( \theta \right) =\log
\phi _{n}\left( \theta \right) =n\psi \left( \theta /n^{1/2}\right)
$.

As we indicated in the introduction, our goal is to extend the
standard Edgeworth expansion given by approximation (\ref{Ed1}),
which is obtained from a formal expansion of the Fourier transform
of $S_{n}/n^{1/2}$. In particular, since $\psi ^{\prime }\left(
0\right) =0$ we have, for fixed $\theta\in\mathbb{R}$,
\begin{equation*}
n\psi \left( \theta /n^{1/2}\right) =\sum_{j=2}^m\frac{\kappa_j(i\theta)^j}{j!n^{j/2-1}}+o\left( \frac{\theta ^{m}}{n^{m/2-1}}\right) ,
\end{equation*}%
if $m<\min \left( \beta ,\gamma \right) $, where $\kappa_j$'s are the cumulants of $X_k$. Therefore,%
\begin{eqnarray}
E\exp \left( \theta S_{n}/n^{1/2}\right) &=&\exp \left( n\psi \left(
\theta
/n^{1/2}\right) \right)  \notag \\
&=&\exp \left( -\theta ^{2}/2\right) \exp \left(
\sum_{j=3}^{m}\frac{\kappa
_{j}\left( i\theta \right) ^{j}}{j!n^{j/2-1}}+o\left( \frac{1}{%
n^{m/2-1}}\right) \right)  \notag \\
&=&\exp \left( -\theta ^{2}/2\right) \left(1+\sum_{k=1}^\infty\frac{1}{k!}\left(\sum_{j=3}^m\frac{\kappa_j(i\theta)^j}{j!n^{j/2-1}}+o\left(\frac{1}{n^{m/2-1}}\right)\right)^k\right) \notag\\
&=&\exp \left( -\theta ^{2}/2\right) \left( 1+\sum_{q=3}^m\sum_{k=1}^q\frac{1}{k!}\xi_{k,q}(\theta)\frac{1}{n^{q/2-1}}+o\left(\frac{1}{n^{m/2-1}}\right)\right) .  \label{EF}
\end{eqnarray}%
where $\xi_{k,q}(\theta)$ is the coefficient of $1/n^{q/2-1}$ in the expansion of $(\sum_{j=3}^q\kappa_j(i\theta)^j/(j!n^{j/2-1}))^k$. Hence the summation in \eqref{EF} has first term $\kappa
_{3}\left( i\theta \right) ^{3}/n^{1/2}$, second term $(\kappa_3^2(i\theta)^6/72+\kappa_4(i\theta)^4/4!)(1/n)$ and so on.

Recall that the inverse Fourier transform of $\left( -i\theta \right)
^{k}\exp \left( -\theta ^{2}/2\right) $ equals $d^{k}\eta \left(
x\right) /dx^{k}$, where $\eta \left( x\right) =\exp \left(
-x^{2}/2\right)/\sqrt{2\pi} $. In
turn, we have that%
\begin{equation*}
\frac{d^{k}}{dx^{k}}\eta \left( x\right) =\eta \left( x\right)
\left( -1\right) ^{k}H_{k}\left( x\right) ,
\end{equation*}%
where $H_{k}\left( x\right) $ is the $k$-th Hermite polynomial. In
particular, $H_{1}\left( x\right) =x$, $H_{2}\left( x\right) =x^{2}-1$ and $%
H_{3}\left( x\right) =x^{3}-3x$ etc. Therefore, taking inverse Fourier
transforms in both sides of (\ref{EF}) we obtain (see Hall (1992) for a derivation of the distribution function, and Evans and Swartz (2000) and Section 3 for density function)
that the density $f_{S_{n}/n^{1/2}}\left( \cdot \right) $ of
$S_{n}/n^{1/2}$
satisfies%
\begin{equation}
f_{S_{n}/n^{1/2}}\left( x\right) =\eta\left( x\right) +\eta\left(
x\right) \sum_{k=3}^{m}n^{-k/2+1}G_{k}\left( x\right) +o\left(
n^{-m/2+1}\right) , \label{ED1}
\end{equation}%
where $G_{k}\left( \cdot \right) $ is a polynomial of degree $k$
that depends only on the first $k$ cumulants of $X_{j}$. In
particular, we have that $G_{3}=\kappa _{3}H_{3}/3!$, $G_{4}=\kappa
_{3}^{2}H_{6}/72+ \kappa _{4} H_{4}/4!$.

Our goal is to understand the contribution of the error term
$o\left( n^{-m/2+1}\right) $ in expansion (\ref{ED1}) and the
problem is that subsequent $G_{k}$'s (for $k\geq\min \left( \beta
,\gamma \right) $) involve moments (or cumulants) that are not
defined for $X_{j}$. The idea is to write $\psi\left( \theta \right)
=\chi\left( \theta \right) +\xi\left( \theta
\right) $ where $\chi\left( \cdot \right) $ is analytic in a neighborhood of the origin while $\xi\left( \cdot \right) $
is a non-analytic function. Dealing with the analytic component, $%
\chi$, gives rise to the standard Edgeworth expansion as in (\ref{ED1}%
), so the interesting part involves the analysis of $\xi\left( \cdot
\right) $ which, as we shall see, yields Theorem \ref{maintheorem}
and \ref{symmetrictheorem} that extends (\ref{ED1}).
\bigskip

\subsection{Main Theorems}

Throughout this paper we assume the following differentiation
condition on $L_\pm$ (for simplicity we write $L$ below):

\begin{assumption} $L$ is bounded, differentiable and satisfies
\begin{equation}
L'(x)=o\left(\frac{L(x)}{x}\right)\label{differentiability}
\end{equation}
as $x\to\infty$ \label{mainassumption}
\end{assumption}

Assumption \ref{mainassumption} controls the order of the derivative
of $L$. Examples of slowly varying function that satisfy
(\ref{differentiability}) include $\log x$, $(\log x)^2$ etc.

We also introduce the function $\zeta_L(x)$ defined as
\begin{equation}
\zeta_L(x)=\int_1^x\frac{L(u)}{u}du
\label{zetadefinition}\end{equation}

Note that Assumption \ref{mainassumption} implies
\begin{equation}
L(x)=o(\zeta_L(x)) \label{integrability}\end{equation} as
$x\to\infty$. Both Assumption \ref{mainassumption} and this
condition will be used in our proof of the main result.

Next let us define Dawson's integral $D(z)$ as
\begin{equation}
D(z)=e^{-z^2}\int_0^z e^{t^2}dt \label{Dawson}
\end{equation}
and the classical parabolic cylinder function $D_\nu(z)$ with parameter $\nu$%
, where $\text{Re}(\nu)>-1$, as
\begin{equation}
D_\nu(z)=\sqrt{\frac{2}{\pi}}e^{z^2/4}\int_0^\infty e^{-t^2/2}t^\nu\cos(zt-%
\frac{\nu\pi}{2})dt \label{classicalparaboliccylinder}
\end{equation}

We state our main theorem:

\begin{thm}
Assume $(X_k: k\geq1)$ are iid rv's with zero mean, unit variance
and density following (\ref{density}) with $L_+$ and $L_-$
satisfying Assumption \ref{mainassumption}. Also assume $\int_{-\infty}^\infty|\phi(\theta)|d\theta<\infty$. We have, for fixed $x>0$,
\begin{equation*}
f_{S_n/n^{1/2}}(x)=\eta(x)\left[1+\sum_{3\leq j<\min(\beta,\gamma)}\frac{G_j(x)}{n^{j/2-1}}\right]%
+F(x,n)+o(F(x,n))
\end{equation*}
as $n\to\infty$, where $G_j(x), j<\min(\beta,\gamma)$ are the ordinary Edgeworth
coefficients (as discussed in Section 2.1) and $F(x,n)$ is defined according to the following
cases:

\textbf{Case 1:\ }$\beta=\gamma=\alpha$

$$
F(x,n)=\left\{\begin{array}{ll}\frac{1}{\Gamma(\alpha+1)n^{\alpha/2-1}}\Big[\frac{e^{-x^2/2}}{\sqrt{2\pi}}H_\alpha(x)\left(\zeta_{L_+}(n^{1/2}x)+\zeta_{L_-}(n^{1/2}x)\right){}&\\
{}-\frac{1}{\sqrt{2}}\frac{d^\alpha}{dx^\alpha}D(\frac{x}{\sqrt{2}})\left(L_+(n^{1/2}x)-L_-(n^{1/2}x)\right)\Big]&\mbox{\
for even
$\alpha$}\\
&\\
\frac{1}{\Gamma(\alpha+1)n^{\alpha/2-1}}\Big[\frac{e^{-x^2/2}}{\sqrt{2\pi}}H_\alpha(x)\left(\zeta_{L_+}(n^{1/2}x)-\zeta_{L_-}(n^{1/2}x)\right){}&\\
{}-\frac{1}{\sqrt{2}}\frac{d^\alpha}{dx^\alpha}D(\frac{x}{\sqrt{2}})\left(L_+(n^{1/2}x)+L_-(n^{1/2}x)\right)\Big]&\mbox{\
for odd
$\alpha$}\\
&\\
-\sqrt{\frac{\pi}{2}}\frac{e^{-x^2/4}}{\Gamma(\alpha+1)\sin(\alpha\pi)n^{\alpha/2-1}}\big[D_\alpha(x)L_+(n^{1/2}x){}&\\
{}+D_\alpha(-x)L_-(n^{1/2}x)\big]&\mbox{\ for non-integer $\alpha$}
\end{array}\right.
$$

\textbf{Case 2:\ }$\beta<\gamma$

$$
F(x,n)=\left\{\begin{array}{ll}\frac{H_\beta(x)e^{-x^2/2}}{\sqrt{2\pi}\Gamma(\beta+1)}\frac{\zeta_{L_+}(n^{1/2}x)}{n^{\beta/2-1}}
&\mbox{\ for integer
$\beta$}\\
-\sqrt{\frac{\pi}{2}}\frac{D_\beta(x)e^{-x^2/4}}{\Gamma(\beta+1)\sin(\beta\pi)}\frac{L_+(n^{1/2}x)}{n^{\beta/2-1}}&\mbox{\
for non-integer $\beta$}
\end{array}\right.
$$

\textbf{Case 3:\ }$\beta>\gamma$

$$
F(x,n)=\left\{\begin{array}{ll}\frac{H_\gamma(x)e^{-x^2/2}}{\sqrt{2\pi}\Gamma(\gamma+1)}\frac{\zeta_{L_-}(n^{1/2}x)}{n^{\gamma/2-1}}
&\mbox{\ for even
$\gamma$}\\
-\frac{H_\gamma(x)e^{-x^2/2}}{\sqrt{2\pi}\Gamma(\gamma+1)}\frac{\zeta_{L_-}(n^{1/2}x)}{n^{\gamma/2-1}}&\mbox{\
for odd
$\gamma$}\\
-\sqrt{\frac{\pi}{2}}\frac{D_\gamma(-x)e^{-x^2/4}}{\Gamma(\gamma+1)\sin(\gamma\pi)}\frac{L_-(n^{1/2}x)}{n^{\gamma/2-1}}&\mbox{\
for non-integer $\gamma$}
\end{array}\right.
$$
Here $H_k(z)$, $D(z)$ and $D_\nu(z)$ are Hermite polynomial of order
$k$, Dawson's integral and classical parabolic cylinder function
with parameter $\nu$ respectively. \label{maintheorem}
\end{thm}

\bigskip

Note that this theorem states the result for $x>0$. For $x<0$, one merely has to consider $-X_k$ and $-S_n/n^{1/2}$ and the result reduces to the $x>0$ case.

This result can be obtained by first considering a symmetric density
and then splitting the non-symmetric density into odd and even
functions which can be tackled separately using result for the
symmetric case. For the symmetric case, i.e. $\beta=\gamma=\alpha$
and $L\equiv L_+=L_-$, we have a neater representation:

\begin{thm}

Assume $(X_k:k\geq1)$ are iid rv's with unit variance and symmetric
density $f(x)=(1+|x|)^{-(\alpha+1)}L(|x|)$ with $L$ satisfying
Assumption \ref{mainassumption} and $\alpha>2$. Moreover, assume $\int_{-\infty}^\infty|\phi(\theta)|d\theta<\infty$. Then, for fixed $x\in\mathbb{R}$,
\begin{equation*}
f_{S_n/n^{1/2}}(x)=\eta(x)\left[1+\sum_{\substack{ 3\leq j<\alpha  \\ j\text{\ even}}}%
\frac{G_j(x)}{n^{j/2-1}}\right]+\left\{\begin{array}{ll}\frac{L(n^{1/2}x)}{n^{\alpha/2-1}}G_\alpha(x)+o(\frac{L(n^{1/2}x)}{n^{\alpha/2-1}})&\text{for non-even $\alpha$}\\
\frac{\zeta_L(n^{1/2}x)}{n^{\alpha/2-1}}G_\alpha(x)+o(\frac{\zeta_L(n^{1/2}x)}{n^{\alpha/2-1}})&\text{for even
$\alpha$}\end{array}\right.
\end{equation*}
as $n\nearrow\infty$, where $\zeta_L$ is defined in
(\ref{zetadefinition}). Here $G_j(x),j<\alpha$ are the ordinary
Edgeworth coefficients (as discussed in Section 2.1), while $G_\alpha(x)$ is defined as
\begin{equation}
G_\alpha(x)=\left\{%
\begin{array}{ll}
-\frac{\sqrt{2}}{\Gamma(\alpha+1)}\frac{d^\alpha}{dx^\alpha}D(\frac{x}{\sqrt{2}})
&
\text{for odd $\alpha$} \\
-\sqrt{\frac{\pi}{2}}\frac{e^{-x^2/4}}{\Gamma(\alpha+1)\sin(\alpha\pi)}%
[D_\alpha(x)+D_\alpha(-x)] & \text{for non-integer $\alpha$}
\\
\sqrt{\frac{2}{\pi}}\frac{e^{-x^2/2}}{\Gamma(\alpha+1)}H_\alpha(x) & \text{%
for even $\alpha$}%
\end{array}%
\right. \label{coefficient}
\end{equation}
where $D(z)$, $D_\nu(z)$ and $H_k(z)$ are Dawson's integral,
classical parabolic cylinder function with parameter $\nu$ and
Hermite polynomial of order $k$ respectively.
\label{symmetrictheorem}
\end{thm}

\bigskip

As we indicated before, the analysis behind the previous result
involves understanding the behavior of the non-analytic component
$\xi$. In fact, we obtain the following decomposition of $\psi$:

\begin{proposition}
Assume a symmetric density as in Theorem \ref{symmetrictheorem}. We
have
\begin{equation}
\psi(\theta)=\chi(\theta)+\xi(\theta)+o(\xi(\theta)) \label{psi}
\end{equation}
where
\begin{equation}
\chi(\theta)=\sum_{\substack{1\leq j<\alpha\\j\text{\
even}}}\frac{(-1)^{j/2}\kappa_j}{j!}\theta^j \label{chi}
\end{equation} is the ordinary Taylor series expansion up to the
largest moment and
\begin{equation}
\xi(\theta)=\left\{\begin{array}{ll}-\frac{\pi}{\Gamma(\alpha+1)\sin
\frac{\alpha\pi}{2}}\ |\theta|^\alpha
L(\frac{1}{|\theta|})&\mbox{\ for non-even $\alpha$}\\
\frac{2(-1)^{\alpha/2}}{\Gamma(\alpha+1)}\ |\theta|^\alpha
\zeta_L(\frac{1}{|\theta|})&\mbox{\ for even
$\alpha$}\end{array}\right. \label{xi}
\end{equation}
is the non-analytic component. Here $\zeta_L$ is defined in
(\ref{zetadefinition}). \label{symmetriccumulant}
\end{proposition}

\bigskip

With this decomposition, we can perform similar derivation as
Equation (\ref{EF}) to obtain
\begin{eqnarray*}
E\exp \left( \theta S_{n}/n^{1/2}\right) &=&\exp \left( -\theta
^{2}/2\right) \left( 1+\sum_{q=3}^m\sum_{k=1}^q\frac{1}{k!}\xi_{k,q}(\theta)\frac{1}{n^{q/2-1}}+n\xi\left(\frac{\theta}{n^{1/2}}\right)+o\left(n\xi\left(\frac{\theta}{n^{1/2}}\right)\right)\right) .
\end{eqnarray*}%

Dawson's integral and the other special functions then arise as the
inverse Fourier transform of $e^{-\theta^2/2}|\theta|^\alpha$
appearing in the non-analytic term of the above expression.

\bigskip

The rest of the paper provides the necessary details behind the analysis of $%
\xi$ and the Fourier inversion required to obtain Theorem
\ref{symmetrictheorem} and then Theorem \ref{maintheorem}; the proof
of Proposition \ref{symmetriccumulant} follows directly from
Proposition \ref{regularlyvaryingcharacteristic}, which is given in
Section 3. As we indicated in the introduction, in Section 4 we
shall discuss the formal connection between our expansion and that
of Rozovskii (1989) in the context of large deviations. It turns out
that the extra term in our expansion, which comes from the
non-analytic component of the characteristic function, blends
smoothly with the heavy-tail asymptotic, at least in some
restrictive cases. When this holds, there exists a smooth transition
from central limit to large deviations region.

\section{Symmetric Density}

Our goal in this section is to prove Theorem \ref{symmetrictheorem}.
As noted before, we adopt the approach of first expanding the
characteristic function of a symmetric regularly varying density,
and hence its cumulant generating function, followed by an inverse
Fourier transform to get the stated result. In order to obtain the
characteristic function for regularly varying density, it helps to
first consider in detail the nature of the approximation in a
simplified context, namely, that of a symmetric Pareto density.

\subsection{Characteristic Function}

\subsubsection{Pareto Density}

In this section we shall assume that%
\begin{equation}
f\left( x\right) =\frac{a_{f}}{1+|x|^{1+\alpha }},\ x\in (-\infty
,\infty ), \label{paretodensity}
\end{equation}%
where $a_{f}$ is the normalizing constant that makes $f\left(
x\right) $ a well defined density. For example, $a_{f}=\pi/\sqrt{2}$
if $\alpha=3$. With this we have the simplest version of our
expansion:

\begin{proposition}
For rv following (\ref{paretodensity}), the characteristic function
takes the form
\begin{equation}
\phi(\theta) =1+ \sum_{\substack{ 2\leq j<\alpha  \\ j \text{\
even}}}\frac{(-1)^{j/2}m_j}{j!}\theta^j +
\left\{\begin{array}{ll}-\frac{\pi a_f}{\Gamma(\alpha+1)\sin
\frac{\alpha\pi}{2}}\ |\theta|^\alpha+o(|\theta|^\alpha)&\mbox{\ for non-even
$\alpha$}\\-\frac{2(-1)^{\alpha/2}a_f}{\Gamma(\alpha+1)}\
|\theta|^\alpha\log|\theta|+o(|\theta|^\alpha\log|\theta|)&\mbox{\ for even
$\alpha$}\end{array}\right. \label{paretoexpansion}
\end{equation}
as $\theta\to0$, where $m_j$ is the $j$-th moment of the rv.
\label{paretocharacteristic}
\end{proposition}

The proof of this theorem is divided into four cases: $\alpha$ is
odd, even, and non-integer with the integral part being odd and
even. Each case involves finding the asymptotic of an integral as
below. Also note that since the density is symmetric, the
characteristic function is real and symmetric. Thus we can consider
only $\theta>0$ without loss of generality. Moreover, throughout the paper we shall use the notation ``$\sim$" to denote equivalence relation i.e. $f(x)\sim g(x)$ iff $f(x)/g(x)\to1$ as $x\searrow0$ (or $x\nearrow\infty$ depending on context). We also use $\lfloor x\rfloor$ as the integral part of $x$ i.e. $\max\{k\in\mathbb{N}:k\leq x\}$ and $\lceil x\rceil$ as the ceiling of $x$ i.e. $\min\{k\in\mathbb{N}:k\geq x\}$.

\begin{lemma} For $\theta\searrow0$, we have the four asymptotics:
\\

\textbf{Case 1:\ } $\alpha$ is odd

\begin{align}
\int_{-\infty}^\infty\frac{x^{\alpha-1}(e^{i\theta
x}-1)}{1+|x|^{1+\alpha}}dx &\sim -\pi\theta \label{oddasymptotic}
\end{align}

\textbf{Case 2:\ } $\alpha$ is even

\begin{align}
\int_{-\infty}^\infty\frac{x^{\alpha-1}e^{i\theta
x}}{1+|x|^{1+\alpha}}dx&\sim -2i\theta\log\theta
\label{evenasymptotic}
\end{align}

For non-integer $\alpha$, denote $q=\llcorner \alpha\lrcorner$. Then
\\

\textbf{Case 3:\ } $q$ is even

\begin{equation*}
\int_{-\infty}^\infty\frac{x^q(e^{i\theta x}-1)}{1+|x|^{1+\alpha}}dx
\sim
2\theta^{\alpha-q}\Gamma(-\alpha+q)\cos\left(\frac{-\alpha+q}{2}\pi\right)
\label{nonintegerevenasymptotic}
\end{equation*}

\textbf{Case 4:\ } $q$ is odd

\begin{equation*}
\int_{-\infty}^\infty\frac{x^qe^{i\theta x}}{1+|x|^{1+\alpha}}dx
\sim
2i\theta^{\alpha-q}\Gamma(-\alpha+q)\sin\left(\frac{-\alpha+q}{2}\pi\right)
\label{nonintegeroddasymptotic}
\end{equation*}

\label{paretolemma}
\end{lemma}

\bigskip

\begin{proof} We prove the lemma case by case:

\begin{itemize}
\item
\textbf{Case 1:\ }$\alpha$ is odd

Since $\alpha$ is odd, the integrand is a symmetric function, and we
can write
\begin{align*}
\int_{-\infty}^\infty\frac{x^{\alpha-1}(e^{i\theta x}-1)}{1+|x|^{1+\alpha}}%
dx&=2\int_0^\infty\frac{x^{\alpha-1}}{1+x^{1+\alpha}}(\cos(\theta
x)-1)dx\\
&=2\theta\int_0^\infty
\frac{u^{\alpha-1}}{\theta^{1+\alpha}+u^{1+\alpha}}(\cos u -1)du
\end{align*}
by a change of variable $u=\theta x$ in the second equality. Since
\begin{equation*}
\left|\frac{u^{\alpha-1}}{\theta^{1+\alpha}+u^{1+\alpha}}(\cos u
-1)\right| \leq \frac{1}{u^2}|\cos u -1|
\end{equation*}
which is integrable, by dominated convergence theorem we get
\begin{equation*}
\int_0^\infty
\frac{u^{\alpha-1}}{\theta^{1+\alpha}+u^{1+\alpha}}(\cos u
-1)du\to\int_0^\infty \frac{\cos u-1}{u^2}du=-\frac{\pi}{2}
\end{equation*}
as $\theta\to 0$. Hence Case 1 is proved.

\item
\textbf{Case 2:\ }$\alpha$ is even

Note that the integrand is an odd function, and
\begin{align*}
\int_{-\infty}^\infty\frac{x^{\alpha-1}e^{i\theta x}}{1+|x|^{1+\alpha}}dx&=2i\int_0^\infty%
\frac{x^{\alpha-1}\sin(\theta x)}{1+x^{1+\alpha}}dx\\
&=2i\theta\int_0^\infty%
\frac{u^{\alpha-1}\sin u}{\theta^{1+\alpha}+u^{1+\alpha}}du
\end{align*}
again by a change of variable $u=\theta x$. Note that as $\theta\to0$, the integral blows up at $u=0$. So for small $\epsilon$%
, we write
\begin{align*}
\int_0^\infty\frac{u^{\alpha-1}\sin u}{\theta^{1+\alpha}+u^{1+\alpha}}%
du&=\int_0^\epsilon\frac{u^{\alpha-1}\sin
u}{\theta^{1+\alpha}+u^{1+\alpha}}du+\int_\epsilon^\infty\frac{u^{\alpha-1}\sin
u}{\theta^{1+\alpha}+u^{1+\alpha}}du
\end{align*}
Now
$$\int_\epsilon^\infty\frac{u^{\alpha-1}\sin
u}{\theta^{1+\alpha}+u^{1+\alpha}}du\leq\int_\epsilon^\infty\frac{1}{u^2}du<\infty$$
On the other hand,
$$\int_0^\epsilon\frac{u^{\alpha-1}\sin
u}{\theta^{1+\alpha}+u^{1+\alpha}}du=\int_0^\epsilon\frac{u^\alpha}{\theta^{1+\alpha}+u^{1+\alpha}}du+R(\theta)$$
where
$$|R(\theta)|\leq C\int_0^\epsilon\frac{u^{2+\alpha}}{\theta^{1+\alpha}+u^{1+\alpha}}du\leq C\int_0^\epsilon udu=O(1)$$
for positive constants $C$ and
$$\int_0^\epsilon\frac{u^\alpha}{\theta^{1+\alpha}+u^{1+\alpha}}du
=\frac{1}{1+\alpha}\log(\theta^{1+\alpha}+u^{1+\alpha})\ \Big|_0^\epsilon\sim-\log\theta$$
Hence Case 2 is proved.

\item
\textbf{Case 3 and 4:\ }$\alpha$ is non-integer with $q$ even and
odd

The proof is similar to the last two cases. For even $q$, we write
\begin{align*}
\int_{-\infty}^\infty\frac{x^q(e^{i\theta
x}-1)}{1+|x|^{1+\alpha}}dx&=2\int_0^\infty
\frac{x^q}{1+x^{1+\alpha}}(\cos(\theta x)-1)dx \\
&=2\theta^{\alpha-q}\int_0^\infty
\frac{u^q}{\theta^{1+\alpha}+u^{1+\alpha}}(\cos
u-1)du \\
&\sim 2\theta^{\alpha-q}\int_0^\infty\frac{\cos
u-1}{u^{1+\alpha-q}}du
\end{align*}
and for odd $q$, we write
\begin{align*}
\int_{-\infty}^\infty\frac{x^qe^{i\theta x}}{1+|x|^{1+\alpha}}dx&=2i\int_0^\infty \frac{%
x^q\sin(\theta x)}{1+x^{1+\alpha}}dx \\
&=2i\theta^{\alpha-q}\int_0^\infty\frac{u^q\sin u}{\theta^{1+\alpha}+u^{1+\alpha}}du \\
&\sim2i\theta^{\alpha-q}\int_0^\infty\frac{\sin u}{u^{1+\alpha-q}}du
\end{align*}
by using a change of variable $u=\theta x$ and dominated convergence
theorem. We use the fact that
\begin{equation*}
\int_0^\infty u^{z-1}(e^{iu}-1)du=\Gamma(z)e^\frac{iz\pi}{2} \text{\
\ for\ } -1<\text{Re}(z)<0
\end{equation*}
to get
\begin{equation*}
\int_0^\infty\frac{\cos u-1}{u^{1+\alpha-q}}du=\Gamma(-\alpha+q)\cos(\frac{%
-\alpha+q}{2}\pi)
\end{equation*}
and
\begin{equation*}
\int_0^\infty\frac{\sin u}{u^{1+\alpha-q}}du=\Gamma(-\alpha+q)\sin(\frac{%
-\alpha+q}{2}\pi)
\end{equation*}
The lemma is then immediate.
\end{itemize}
\end{proof}

\bigskip

We are now ready to prove Proposition \ref{paretocharacteristic}.
The proof basically reduces (\ref{paretoexpansion}) into expressions
in terms of the above four integrals through multiple use of L'
Hospital's rule.

\begin{proof}[Proof of Proposition \ref{paretocharacteristic}] It suffices to consider $\theta>0$. Again we
consider case by case:

\begin{itemize}
\item
\textbf{Case 1:\ }$\alpha$ is odd

Applying L'Hospital's rule successively for $\alpha-1$ times, we
have
\begin{eqnarray*}
& &\lim_{\theta\to0}\frac{\int_{-\infty}^\infty \frac{a_fe^{i\theta x}}{1+|x|^{1+\alpha}}%
dx-1-\sum_{\substack{ 1\leq j<\alpha  \\ j \text{\ even}}}\frac{(-1)^{j/2}m_j%
}{j!}\theta^j}{\theta^\alpha} \\
&=& \lim_{\theta\to0}\frac{\int_{-\infty}^\infty\frac{a_fixe^{i\theta x}}{1+|x|^{1+\alpha}}%
dx-\sum_{\substack{ 1\leq j <\alpha  \\ j \text{\ even}}}\frac{(-1)^{j/2}m_j%
}{(j-1)!}\theta^{j-1}}{\alpha \theta^{\alpha-1}}\\
&=&\lim_{\theta\to0}\frac{\int_{-\infty}^\infty\frac{a_fi^2x^2e^{i\theta x}}{1+|x|^{1+\alpha}}%
dx-\sum_{\substack{ 2\leq j <\alpha  \\ j \text{\ even}}}\frac{(-1)^{j/2}m_j%
}{(j-2)!}\theta^{j-2}}{\alpha(\alpha-1) \theta^{\alpha-2}}\\
&&\vdots\\
&=&\lim_{\theta\to0}\frac{\int_{-\infty}^\infty\frac{a_fi^{\alpha-1}x^{\alpha-1}e^{i\theta x}}{%
1+|x|^{1+\alpha}}dx-(-1)^{(\alpha-1)/2}m_{\alpha-1}}{\alpha!\theta} \\
&=&\frac{(-1)^{(\alpha-1)/2}a_f}{\alpha!}\lim_{\theta\to0}\frac{\int_{-\infty}^\infty\frac{%
x^{\alpha-1}(e^{i\theta x}-1)}{1+|x|^{1+\alpha}}dx}{\theta} \\
&=&\frac{(-1)^{(\alpha+1)/2}\pi a_f}{\alpha!}
\end{eqnarray*}
where the last step follows from (\ref{oddasymptotic}). Note that
the integral in each step above is finite since the moments exist up
to order $\alpha-1$.

The result for odd $\alpha$ thus follows.

\item
\textbf{Case 2:\ }$\alpha$ is even

Using the fact that
\begin{equation*}
\frac{d^m}{dk^m}k^\alpha\log
k=\alpha(\alpha-1)\cdots(\alpha-m+1)k^{\alpha-m}\log k +
\gamma_{\alpha-m}k^{\alpha-m}
\end{equation*}
for $m=1,\ldots,\alpha-1$, where $\gamma_{\alpha-m}$ are constants,
and applying L' Hospital's rule for $\alpha-1$ times, we get
\begin{eqnarray*}
& &\lim_{\theta\to0}\frac{\int_{-\infty}^\infty \frac{a_fe^{i\theta x}}{1+|x|^{1+\alpha}}%
dx-1-\sum_{\substack{ 1\leq j<\alpha  \\ j \text{\ even}}}\frac{(-1)^{j/2}m_j%
}{j!}\theta^j}{\theta^\alpha\log\theta} \\
&=&\lim_{\theta\to0}\frac{\int_{-\infty}^\infty\frac{a_fi^{\alpha-1}x^{\alpha-1}e^{i\theta x}}{%
1+|x|^{1+\alpha}}dx}{\alpha!\theta\log\theta+\gamma_1\theta} \\
&=&\frac{a_fi^{\alpha-1}}{\alpha!}\lim_{\theta\to0}\frac{\int_{-\infty}^\infty\frac{%
x^{\alpha-1}e^{i\theta x}}{1+|x|^{1+\alpha}}dx}{\theta\log\theta} \\
&=&-\frac{2(-1)^{\alpha/2}a_f}{\alpha!}
\end{eqnarray*}
where the last equality follows from (\ref{evenasymptotic}).

The result for even $\alpha$ thus follows.

\item
\textbf{Case 3 and 4:\ }$\alpha$ is non-integer with $q$ even and
odd

Following the same line of proof as the previous two cases, we get,
for even $q$,
\begin{eqnarray*}
& &\lim_{\theta\to0}\frac{\int_{-\infty}^\infty \frac{a_fe^{i\theta x}}{1+|x|^{1+\alpha}}%
dx-1-\sum_{\substack{ 1\leq j<\alpha  \\ j \text{\ even}}}\frac{(-1)^{j/2}m_j%
}{j!}\theta^j}{\theta^\alpha} \\
&=&\lim_{\theta\to0}\frac{i^qa_f\int_{-\infty}^\infty\frac{x^q(e^{i\theta x}-1)}{%
1+|x|^{1+\alpha}}dx}{\alpha(\alpha-1)\cdots(\alpha-q+1)\theta^{\alpha-q}} \\
&=&\frac{(-1)^{q/2}a_f\Gamma(\alpha-q+1)}{\Gamma(\alpha+1)}\cdot
2\Gamma(-\alpha+q)\cos(\frac{-\alpha+q}{2}\pi) \\
&=&\frac{2(-1)^{q/2}\pi a_f\cos(\frac{-\alpha+q}{2}\pi)}{\Gamma(\alpha+1)%
\sin((-\alpha+q)\pi)} \\
&=&-\frac{(-1)^{q/2}\pi a_f}{\Gamma(\alpha+1)\sin(\frac{\alpha-q}{2}\pi)} \\
&=&-\frac{(-1)^{q/2}\pi a_f}{\Gamma(\alpha+1)\sin\frac{\alpha\pi}{2}\cos\frac{%
q\pi}{2}} \\
&=&-\frac{\pi a_f}{\Gamma(\alpha+1)\sin\frac{\alpha\pi}{2}}
\end{eqnarray*}

and for odd $q$,

\begin{eqnarray*}
& &\lim_{\theta\to0}\frac{\int_{-\infty}^\infty \frac{a_fe^{i\theta x}}{1+|x|^{1+\alpha}}%
dx-1-\sum_{\substack{ 1\leq j<\alpha  \\ j \text{\ even}}}\frac{(-1)^{j/2}m_j%
}{j!}\theta^j}{\theta^\alpha} \\
&=&\lim_{\theta\to0}\frac{i^qa_f\int_{-\infty}^\infty\frac{x^qe^{i\theta x}}{1+|x|^{1+\alpha}}dx%
}{\alpha(\alpha-1)\cdots(\alpha-q+1)\theta^{\alpha-q}} \\
&=&\frac{(-1)^{(q+1)/2}a_f\Gamma(\alpha-q+1)}{\Gamma(\alpha+1)}\cdot
2\Gamma(-\alpha+q)\sin(\frac{-\alpha+q}{2}\pi) \\
&=&\frac{2(-1)^{(q+1)/2}\pi a_f\sin(\frac{-\alpha+q}{2}\pi)}{%
\Gamma(\alpha+1)\sin((-\alpha+q)\pi)} \\
&=&\frac{(-1)^{(q+1)/2}\pi a_f}{\Gamma(\alpha+1)\cos(\frac{\alpha-q}{2}\pi)} \\
&=&\frac{(-1)^{(q+1)/2}\pi a_f}{\Gamma(\alpha+1)\sin\frac{\alpha\pi}{2}\sin%
\frac{q\pi}{2}} \\
&=&-\frac{\pi a_f}{\Gamma(\alpha+1)\sin\frac{\alpha\pi}{2}}
\end{eqnarray*}
by using the reflection property of gamma function and also the
double-angle and sum-of-angle trigonometric identities.

We have proven the proposition.
\end{itemize}
\end{proof}

\subsubsection{Regularly Varying Density}

We are going to extend the result in Proposition
\ref{paretocharacteristic} to density with regularly varying tail in
the form
\begin{equation}
f\left( x\right) =\frac{L(|x|)}{1+|x|^{1+\alpha }},\ x\in (-\infty
,\infty ) \label{regularlyvaryingdensity}
\end{equation}%
where $L$ is a slowly regularly varying function and satisfies the
differentiation condition in Assumption \ref{mainassumption}.

\begin{proposition}
For rv following (\ref{regularlyvaryingdensity}) with $L$ satisfying
Assumption \ref{mainassumption}, the characteristic function takes
the form
\begin{equation}
\phi(\theta) = 1+\sum_{\substack{ 2\leq j<\alpha  \\ j \text{\
even}}}\frac{(-1)^{j/2}m_j}{j!}\theta^j +
\left\{\begin{array}{ll}-\frac{\pi }{\Gamma(\alpha+1)\sin
\frac{\alpha\pi}{2}}\ |\theta|^\alpha
L\left(\frac{1}{|\theta|}\right)+o\left(|\theta|^\alpha L\left(\frac{1}{|\theta|}\right)\right)&\mbox{\ for non-even
$\alpha$}\\\frac{2(-1)^{\alpha/2}}{\Gamma(\alpha+1)}\
|\theta|^\alpha\zeta_L\left(\frac{1}{|\theta|}\right)+o\left(|\theta|^\alpha \zeta_L\left(\frac{1}{|\theta|}\right)\right)&\mbox{\ for
even $\alpha$}\end{array}\right. \label{regularlyvaryingexpansion}
\end{equation}
as $\theta\to0$, where $m_j$ is the $j$-th moment of the rv and $\zeta_L$ is defined
in (\ref{zetadefinition}). \label{regularlyvaryingcharacteristic}
\end{proposition}

\bigskip

To prove Proposition \ref{regularlyvaryingcharacteristic}, we need
the following property of slowly varying function:

\begin{lemma}
Let $L$ be a slowly varying function. Then for any $\rho_1,\rho_2>0$, we can pick an $%
\eta>0$ such that for $0<\theta<\eta$,
\begin{equation*}
\frac{L(\frac{x}{\theta})}{L(\frac{1}{\theta})}\leq\left\{%
\begin{array}{ll}
Cx^{\rho_1} & \text{for\ \ }x>1 \\
Cx^{-\rho_2} & \text{for\ \ }x\leq1%
\end{array}%
\right.
\end{equation*}
where $C$ is a positive constant. \label{regularlyvaryinglemma}
\end{lemma}

\bigskip

\begin{proof} Karamata's representation (see for example, Resnick, 1987) states that for
a slowly varying function $L:\mathbb{R}_+\to\mathbb{R}%
_+$, we have
\begin{equation*}
L(x)=c(x)\exp\left\{\int_1^x t^{-1}\epsilon(t)dt\right\}
\end{equation*}
where $c:\mathbb{R}_+\to\mathbb{R}_+$ and
$\epsilon:\mathbb{R}\to\mathbb{R}$ such that
\begin{align*}
\lim_{x\to\infty} c(x)&=c\in(0,\infty) \\
\lim_{t\to\infty}\epsilon(t)&=0
\end{align*}
So we have
\begin{equation*}
\frac{L(\frac{x}{\theta})}{L(\frac{1}{\theta})}=\frac{c(\frac{x}{\theta})}{c(\frac{1}{\theta})}%
\exp\left\{\int_\frac{1}{\theta}^\frac{x}{\theta}t^{-1}\epsilon(t)dt\right\}
\end{equation*}
Since $L$ is bounded, we can choose $c(x)$ bounded from above and $c(1/\theta)>c_1>0$ for
$\theta$ small enough. So
\begin{equation*}
\frac{c(\frac{x}{\theta})}{c(\frac{1}{\theta})}\leq C
\end{equation*}
for some $C>0$ for small $\theta$. Also by a change of variable
$t=s/\theta$, we have
\begin{equation*}
\exp\left\{\int_\frac{1}{\theta}^\frac{x}{\theta}t^{-1}\epsilon(t)dt\right\}=\exp\left%
\{\int_1^xs^{-1}\epsilon(\frac{s}{\theta})ds\right\}
\end{equation*}
Suppose $x>1$. For any given $\rho_1$, if we pick $\theta$ small enough, we have $%
\epsilon(\frac{s}{\theta})<\rho_1$ for all $s>1$, and
\begin{equation*}
\exp\left\{\int_1^xs^{-1}\epsilon(\frac{s}{\theta})ds\right\}\leq\exp\left\{%
\int_1^xs^{-1}\rho_1 ds\right\}=x^{\rho_1}
\end{equation*}
The result for $x>1$ follows. Similarly, for $x\leq1$, if we pick
$\theta$ small enough we have $\epsilon(\frac{s}{\theta})>-\rho_2$,
and
\begin{equation*}
\exp\left\{\int_1^xs^{-1}\epsilon(\frac{s}{\theta})ds\right\}=\exp\left%
\{-\int_x^1s^{-1}\epsilon(\frac{s}{\theta})ds\right\}\leq\exp\left\{\int_x^1s^{-1}\rho_2ds\right\}\leq
x^{-\rho_2}
\end{equation*}
Thus the lemma is proved.
\end{proof}

\bigskip

We now state the correspondence of Lemma \ref{paretolemma} for
regularly varying density:

\begin{lemma} For $\theta\searrow0$, we have the four asymptotics:
\\

\textbf{Case 1:\ }$\alpha$ is odd
\begin{equation*}
\int_{-\infty}^\infty\frac{x^{\alpha-1}L(|x|)(e^{i\theta
x}-1)}{1+|x|^{1+\alpha}}dx \sim -\pi \theta
L\left(\frac{1}{\theta}\right)
\end{equation*}

\textbf{Case 2:\ }$\alpha$ is even
\begin{equation*}
\int_{-\infty}^\infty\frac{x^{\alpha-1}L(|x|)e^{i\theta
x}}{1+|x|^{1+\alpha}}dx\sim2i\theta\zeta_L\left(\frac{1}{\theta}\right)
\end{equation*}

where $\zeta_L$ is defined in (\ref{zetadefinition}). For
non-integer $\alpha$, denote $q=\llcorner \alpha\lrcorner$ as the
integral part of $\alpha$. Then
\\

\textbf{Case 3:\ } $q$ is even

\begin{equation*}
\int_{-\infty}^\infty\frac{x^qL(|x|)(e^{i\theta
x}-1)}{1+|x|^{1+\alpha}}dx \sim
2\theta^{\alpha-q}L\left(\frac{1}{\theta}\right)\Gamma(-\alpha+q)\cos\left(\frac{-\alpha+q}{2}\pi\right)
\end{equation*}

\textbf{Case 4:\ } $q$ is odd

\begin{equation*}
\int_{-\infty}^\infty\frac{x^qL(|x|)e^{i\theta
x}}{1+|x|^{1+\alpha}}dx \sim
2i\theta^{\alpha-q}L\left(\frac{1}{\theta}\right)\Gamma(-\alpha+q)\sin\left(\frac{-\alpha+q}{2}\pi\right)
\end{equation*}

\end{lemma}

\bigskip

\begin{proof}
Again we prove the lemma case by case.

\begin{itemize}
\item
\textbf{Case 1: }$\alpha$ is odd

We follow similar line of proof as in lemma \ref{paretolemma}. Using
a change of variables $u=\theta x$ we obtain
\begin{align*}
\int_{-\infty}^\infty\frac{x^{\alpha-1}L(|x|)(e^{i\theta x}-1)}{1+|x|^{1+\alpha}}%
dx&=2\theta\int_0^\infty \frac{u^{\alpha-1}L(\frac{u}{\theta})}{\theta^{1+\alpha}+u^{1+%
\alpha}}(\cos u -1)du \\
&=2\theta L(\frac{1}{\theta})\int_0^\infty\frac{u^{\alpha-1}}{\theta^{1+\alpha}+u^{1+\alpha}}%
\cdot\frac{L(\frac{u}{\theta})}{L(\frac{1}{\theta})}\cdot(\cos
u-1)du
\end{align*}
Now using lemma \ref{regularlyvaryinglemma}, for $\theta$ small
enough, we can pick $0<\rho_1,\rho_2<1$ such that
\begin{eqnarray*}
&&\frac{u^{\alpha-1}}{\theta^{1+\alpha}+u^{1+\alpha}}\cdot\frac{L(\frac{u}{\theta})}{L(%
\frac{1}{\theta})}\cdot(\cos u-1) \\
&\leq&\frac{C}{u^{2+\rho_2}}(\cos
u-1)I(u\leq1)+\frac{C}{u^{2-\rho_1}}(\cos u-1)I(u>1)
\end{eqnarray*}
for some $C>0$, for small enough $\theta$. Note that the majorizing function is integrable.
Also note that by the definition of slow variation we have
\begin{equation*}
\frac{L(\frac{u}{\theta})}{L(\frac{1}{\theta})}\to1
\end{equation*}
as $\theta\to0$ for all $u>0$. So by dominated convergence theorem,
\begin{equation*}
\int_0^\infty\frac{u^{\alpha-1}}{\theta^{1+\alpha}+u^{1+\alpha}}\cdot\frac{L(%
\frac{u}{\theta})}{L(\frac{1}{\theta})}\cdot(\cos u-1)du\to\int_0^\infty\frac{1}{u^2}%
(\cos u-1)du
\end{equation*}
and the asymptotic follows.

\item
\textbf{Case 2:\ }$\alpha$ is even

We write
\begin{align*}
\int_{-\infty}^\infty\frac{x^{\alpha-1}L(|x|)e^{i\theta
x}}{1+|x|^{1+\alpha}}dx&=2i\int_0^\infty\frac{x^{\alpha-1}L(x)\sin(\theta
x)}{1+x^{1+\alpha}}dx\\
&=2i\theta\int_0^\infty\frac{u^{\alpha-1}L(\frac{u}{\theta})\sin
u}{\theta^{1+\alpha}+u^{1+\alpha}}du\\
&=2i\theta\left[\int_0^\epsilon\frac{u^{\alpha-1}L(\frac{u}{\theta})\sin
u}{\theta^{1+\alpha}+u^{1+\alpha}}du+\int_\epsilon^\infty\frac{u^{\alpha-1}L(\frac{u}{\theta})\sin
u}{\theta^{1+\alpha}+u^{1+\alpha}}du\right]
\end{align*}
for some small $\epsilon>0$. For the second integral, if we choose $0<\rho_1<1$ in Lemma
\ref{regularlyvaryinglemma}, we have the following bound
\begin{equation}\left|\int_\epsilon^\infty\frac{u^{\alpha-1}L(\frac{u}{\theta})\sin
u}{\theta^{1+\alpha}+u^{1+\alpha}}du\right|\leq
L\left(\frac{1}{\theta}\right)\int_\epsilon^\infty\frac{u^{\alpha-1}Cu^{\rho_1}}{\theta^{1+\alpha}+u^{1+\alpha}}du\leq
M_1L\left(\frac{1}{\theta}\right) \label{tailbound}\end{equation}
for some $C,M_1>0$, for small enough $\theta$.

For the first integral, we write
$$\int_0^\epsilon\frac{u^{\alpha-1}L(\frac{u}{\theta})\sin
u}{\theta^{1+\alpha}+u^{1+\alpha}}du=\int_0^\epsilon\frac{u^\alpha
L(\frac{u}{\theta})}{\theta^{1+\alpha}+u^{1+\alpha}}du+R(\theta)$$
where \begin{equation}|R(\theta)|\leq
M_2L\left(\frac{1}{\theta}\right)\int_0^\epsilon\frac{u^{\alpha+2}u^{-\rho_2}}{\theta^{1+\alpha}+u^{1+\alpha}}du\leq
M_3L\left(\frac{1}{\theta}\right)
\label{secondorderbound}\end{equation} for some $M_2,M_3>0$, $0<\rho_2<2$ and $\theta$ small enough.

Now do a change of variable $x=u/\theta$ back, we get
\begin{equation}\int_0^1\frac{u^\alpha
L(\frac{u}{\theta})}{\theta^{1+\alpha}+u^{1+\alpha}}du=\int_0^\frac{1}{\theta}\frac{x^\alpha
L(x)}{1+x^{1+\alpha}}dx\sim\int_1^\frac{1}{\theta}\frac{L(x)}{x}dx
\label{asymptotic}\end{equation}

Combining (\ref{tailbound}), (\ref{secondorderbound}) and
(\ref{asymptotic}), we obtain the result for even $\alpha$.

\item
\textbf{Case 3 and 4:\ }$\alpha$ is non-integer with $q$ even and
odd

Using the same line of proof as Lemma \ref{paretolemma} and Case 1, and taking $0<\rho_1<\alpha-q$ and $0<\rho_2<2-\alpha+q$, we get the stated result.

\end{itemize}
\end{proof}

The proof of Proposition \ref{regularlyvaryingcharacteristic} is
similar to that of Pareto density, and hence we only outline the proof below. Proposition \ref{symmetriccumulant} then
follows immediately by letting $\psi=\log\phi$ and applying a Taylor
expansion. In doing so one just has to verify that the non-analytic
term in $\phi$ is the only term of that order in the expansion of
$\psi$ to conclude the result.

\begin{proof}[Outline of Proof of Proposition
\ref{regularlyvaryingcharacteristic}] Again it suffices to consider
$\theta>0$, and we use L' Hospital's rule successively. For odd
$\alpha$, recognizing that $L'(\frac{1}{\theta})=o(\theta
L(\frac{1}{\theta}))$ as $\theta\to0$ by Assumption
\ref{mainassumption}, we have
\begin{eqnarray*}
& &\lim_{\theta\to0}\frac{\int_{-\infty}^\infty \frac{L(x)e^{i\theta x}}{1+|x|^{1+\alpha}}%
dx-1-\sum_{\substack{ 1\leq j<\alpha  \\ j \text{\ even}}}\frac{(-1)^{j/2}m_j%
}{j!}\theta^j}{\theta^\alpha L(\frac{1}{\theta})} \\
&=& \lim_{\theta\to0}\frac{\int_{-\infty}^\infty\frac{L(x)ixe^{i\theta x}}{1+|x|^{1+\alpha}}%
dx-\sum_{\substack{ 1\leq j <\alpha  \\ j \text{\ even}}}\frac{(-1)^{j/2}m_j%
}{(j-1)!}\theta^{j-1}}{\alpha \theta^{\alpha-1}L(\frac{1}{\theta})-\theta^{\alpha-2}L'(\frac{1}{\theta})} \\
&=& \lim_{\theta\to0}\frac{\int_{-\infty}^\infty\frac{L(x)ixe^{i\theta x}}{1+|x|^{1+\alpha}}%
dx-\sum_{\substack{ 1\leq j <\alpha  \\ j \text{\ even}}}\frac{(-1)^{j/2}m_j%
}{(j-1)!}\theta^{j-1}}{\alpha \theta^{\alpha-1}L(\frac{1}{\theta})}\\
&&\vdots \\
&=&\frac{(-1)^{(\alpha-1)/2}}{\alpha!}\lim_{\theta\to0}\frac{\int_{-\infty}^\infty\frac{%
L(x)x^{\alpha-1}(e^{i\theta x}-1)}{1+|x|^{1+\alpha}}dx}{\theta L(\frac{1}{\theta})} \\
&=&\frac{(-1)^{(\alpha+1)/2}\pi}{\alpha!}
\end{eqnarray*}

The case for even $\alpha$ can be proved similarly by using
(\ref{integrability}) that states that
$L(\frac{1}{\theta})=o(\zeta_L(\frac{1}{\theta}))$, while the case
for non-integer $\alpha$ is proved similarly with Assumption
\ref{mainassumption}.

\end{proof}

\subsection{Beyond Edgeworth Expansion}

We will now prove the main theorem for the symmetric case, namely
Theorem \ref{symmetrictheorem}. First of all we need to evaluate the
following integrals, which corresponds to a Fourier inversion of the
non-analytic component of the expansion studied in Section 3.1:

\begin{lemma}
\begin{equation*}
\int_{-\infty}^\infty e^{-i\theta x-\theta^2/2}|\theta|^\alpha
d\theta=\left\{\begin{array}{ll}2\sqrt{2}(-1)^{(\alpha-1)/2}\frac{d^\alpha}{dx^\alpha}D(\frac{x}{\sqrt{2}})&\mbox{\
for odd $\alpha$}\\
\sqrt{\frac{\pi}{2}}\sec\frac{\alpha\pi}{2}e^{-x^2/4}[D_\alpha(x)+D_\alpha(-x)]&\mbox{\
for non-integer $\alpha$}\\
(-1)^{\alpha/2}\sqrt{2\pi}e^{-x^2/2}H_\alpha(x)&\mbox{\ for even
$\alpha$}\end{array}\right.
\end{equation*}
where $D(z)$, $D_\nu(z)$ and $H_k(z)$ are Dawson's integral,
classical parabolic cylinder function with parameter $\nu$ and
Hermite polynomial of order $k$ respectively.
\label{Edgeworthlemma}\end{lemma}

\bigskip

\begin{proof} Consider the following
cases:
\begin{itemize}
\item
\textbf{Case 1:\ }$\alpha$ is odd

Observe that
\begin{align*}
\int_{-\infty}^\infty e^{-i\theta x-\theta^2/2}|\theta|^\alpha
d\theta&=2\int_0^\infty\cos(\theta
x)e^{-\theta^2/2}\theta^\alpha d\theta\\
&=2(-1)^{(\alpha-1)/2}\frac{%
d^\alpha}{dx^\alpha}\int_0^\infty\sin(\theta x)e^{-\theta^2/2}d\theta \\
&=2(-1)^{(\alpha-1)/2}\frac{d^\alpha}{dx^\alpha}\text{Im}\left[\int_0^\infty
e^{i\theta x-\theta^2/2}d\theta\right] \\
&=2(-1)^{(\alpha-1)/2}\frac{d^\alpha}{dx^\alpha}\text{Im}\left[%
e^{-x^2/2}\int_{-ix}^{\infty-ix}e^{-q^2/2}dq\right]
\end{align*}
by letting $q=\theta-ix$. We evaluate the integral
\begin{equation*}
\int_{-ix}^{\infty-ix}e^{-q^2/2}dq
\end{equation*}
by forming a closed anti-clockwise rectangular contour from $-ix$ to $R-ix$ to $R$ to 0 to $%
-ix$. The contour integral is zero by analyticity. Note that
\begin{equation*}
\int_0^\infty e^{-q^2/2}dq=\sqrt{\frac{\pi}{2}}
\end{equation*}
The integral from $R-ix$ to $R$ tends to zero as $R\to\infty$. By a
change of variable $t=q/i\sqrt{2}$, the integral along the segment
from $0$ to $-ix$ is
\begin{align*}
\int_0^{-ix}e^{-q^2/2}dq&=i\sqrt{2}\int_0^{-x/\sqrt{2}}e^{t^2}dt \\
&=-i\sqrt{2}e^{x^2/2}D(\frac{x}{\sqrt{2}})
\end{align*}
where $D(z)$ is the Dawson's integral. Therefore we have
\begin{equation*}
\int_{-ix}^{\infty-ix}e^{-q^2/2}dq=\sqrt{\frac{\pi}{2}}+i\sqrt{2}e^{x^2/2}D(%
\frac{x}{\sqrt{2}})
\end{equation*}
and the result follows.

\item
\textbf{Case 2:\ }$\alpha$ is non-integer

Again we have $$\int_{-\infty}^\infty e^{-i\theta
x-\theta^2/2}|\theta|^\alpha d\theta=2\int_0^\infty\cos(\theta
x)e^{-\theta^2/2}\theta^\alpha d\theta$$ Using sum-of-angle formula
and (\ref{classicalparaboliccylinder}), we can write the classical
parabolic cylinder function as
\begin{equation*}
D_\alpha(x)=\sqrt{\frac{2}{\pi}}\cos\frac{\alpha\pi}{2}e^{x^2/4}\int_0^%
\infty e^{-t^2/2}t^\alpha\cos(tx)dt+\sqrt{\frac{2}{\pi}}\sin\frac{\alpha\pi}{%
2}e^{x^2/4}\int_0^\infty e^{-t^2/2}t^\alpha\sin(tx)dt
\end{equation*}
Since the first term is even while the second term is odd in $x$, we
have
\begin{equation*}
D_\alpha(x)+D_\alpha(-x)=2\sqrt{\frac{2}{\pi}}\cos\frac{\alpha\pi}{2}%
e^{x^2/4}\int_0^\infty e^{-t^2/2}t^\alpha\cos(tx)dt
\end{equation*}
The result follows.

\item
\textbf{Case 3:\ }$\alpha$ is even

Note that the integral now becomes analytic and the result is
obvious by the discussion in Section 2.1.

\end{itemize}
\end{proof}

We are now ready to prove Theorem \ref{symmetrictheorem}:

\begin{proof}[Proof of Theorem \ref{symmetrictheorem}] We follow the
proof outline of Evans and Swartz (2000). First note that we can express the
characteristic function
$$\psi(\theta)=\chi(\theta)+\xi(\theta)+R(\theta)$$
where $\chi$ and $\xi$ are as defined in (\ref{chi}) and (\ref{xi}),
and $R$ is the error term. Also let
$$\lambda(\theta)=\chi(\theta)+\frac{\theta^2}{2}+\xi(\theta)$$ be the expansion of characteristic function but trimming the first term and the error term.

Let
$$w_n(\theta)=\sum_{k=0}^{\lfloor\alpha/2-1\rfloor}\frac{(n\lambda(\theta/\sqrt{n}))^k}{k!}$$
and $q_n(\theta)$ be the expansion of $w_n(\theta)$ up to order
$$\left\{\begin{array}{ll}1/n^{\alpha/2-1}L(\frac{n^{1/2}}{|\theta|})&\mbox{for non-even $\alpha$}\\
1/n^{\alpha/2-1}\zeta_L(\frac{n^{1/2}}{|\theta|})&\mbox{for even $\alpha$}
\end{array}\right.$$
In other words,
$$q_n(\theta)=1+\sum_{4\leq j<\alpha}\sum_{k=1}^j\frac{1}{k!}\xi_{k,j}(\theta)\frac{1}{n^{j/2-1}}+\left\{\begin{array}{ll}-\frac{\pi}{\Gamma(\alpha+1)\sin
\frac{\alpha\pi}{2}}\ \frac{|\theta|^\alpha
L(n^{1/2}/|\theta|)}{n^{\alpha/2-1}}&\mbox{\ for non-even $\alpha$}\\
\frac{2(-1)^{\alpha/2}}{\Gamma(\alpha+1)}\ \frac{|\theta|^\alpha
\zeta_L(n^{1/2}/|\theta|)}{n^{\alpha/2-1}}&\mbox{\ for even
$\alpha$}\end{array}\right.$$
where $\xi_{k,j}(\theta)$'s are defined as in \eqref{EF}. Also let
$r_n(\theta)=w_n(\theta)-q_n(\theta)$ be the error of this trimming.
It is easy to verify through direct expansion that $r_n(\theta)$ is a sum consisting of terms written as $P_k(\theta,n)n^{-k}$, where $k\geq\alpha/2-1$ and $P_k(\theta,n)$ is a sum of terms in the form
$$c|\theta|^b\times\left\{\begin{array}{ll}(\frac{L(n^{1/2})}{|\theta|})^d&\mbox{for non-even $\alpha$}\\
(\frac{\zeta_L(n^{1/2})}{|\theta|})^d&\mbox{for even $\alpha$}
\end{array}\right.$$
for $c\in\mathbb{R}$, $b>0$ and $d\in\mathbb{N}$ (and note that for $k=\alpha/2-1$ the terms in $P_k(\theta,n)$ must have $d=0$).

We first show that
$$\int_{-\infty}^\infty e^{-i\theta
x-\theta^2/2}q_n(\theta)\frac{d\theta}{2\pi}$$ is a good
approximation to $f_{S_n/\sqrt{n}}(x)=\int_{-\infty}^\infty
e^{-i\theta x+n\psi(\theta/\sqrt{n})}d\theta/2\pi$. After that we
will verify using Lemma \ref{Edgeworthlemma} that this approximation
is the one shown in Theorem \ref{symmetrictheorem}. Consider
\begin{eqnarray}
&&\left|\int_{-\infty}^\infty\left[e^{-i\theta
x+n\psi(\theta/\sqrt{n})}-e^{-i\theta
x-\theta^2/2}q_n(\theta)\right]\frac{d\theta}{2\pi}\right| \notag\\
&\leq&\int_{-\delta\sqrt{n}}^{\delta\sqrt{n}}|e^{n\psi(\theta/\sqrt{n})}-e^{-\theta^2/2}q_n(\theta)|\frac{d\theta}{2\pi}+{} \notag\\
&&{}\int_{\theta\notin(-\delta\sqrt{n},\delta\sqrt{n})}
e^{-\theta^2/2}|q_n(\theta)|\frac{d\theta}{2\pi}+\int_{\theta\notin(-\delta\sqrt{n},\delta\sqrt{n})}\left|\phi\left(\frac{\theta}{\sqrt{n}}\right)\right|^n\frac{d\theta}{2\pi} \label{intermediate1}
\end{eqnarray}
for any $\delta$. Later we are going to choose the value of $\delta$
to bound our error.

We first consider the third integral in \eqref{intermediate1}. By non-lattice assumption,
given any $\delta$, we can find $\tau_\delta<1$ such that
$|\phi(y)|\leq\tau_\delta$ for any $y>\delta$. Hence
$$\int_{\theta\notin(-\delta\sqrt{n},\delta\sqrt{n})}\left|\phi\left(\frac{\theta}{\sqrt{n}}\right)\right|^n\frac{d\theta}{2\pi}\leq\tau_\delta^{n-1}\int_{\theta\notin(-\delta\sqrt{n},\delta\sqrt{n})}\left|\phi\left(\frac{\theta}{\sqrt{n}}\right)\right|\frac{d\theta}{2\pi}=o(n^{-p})$$
for any $p>0$. The last equation follows from our assumption that
$\phi$ is integrable.

Now consider the second integral in \eqref{intermediate1}. Note that for any $\beta\geq0$,
\begin{equation}
\int_{\theta\notin(-\delta\sqrt{n},\delta\sqrt{n})}e^{-\theta^2/2}\theta^\beta
d\theta =o(n^{-p}) \label{intermediate5}
\end{equation}
for any $p>0$. Also, for any $\alpha>2$ and
slowly varying function $U$,
\begin{equation}
\int_{\theta\notin(-\delta\sqrt{n},\delta\sqrt{n})}e^{-\theta^2/2}|\theta|^\alpha
U\left(\frac{n^{1/2}}{|\theta|}\right)d\theta\leq C_q
\int_{\theta\notin(-\delta\sqrt{n},\delta\sqrt{n})}e^{-\theta^2/2}|\theta|^{\alpha-q}n^{q/2}d\theta=o(n^{-p}) \label{intermediate6}
\end{equation}
for any $q>0$ and $C_q$ a positive constant depending on $q$. Again
$p>0$ is arbitrary. By plugging in the expansion of $q_n$, together
with observations \eqref{intermediate5} and \eqref{intermediate6}, we conclude that the second
integral in \eqref{intermediate1} is $o(n^{-p})$ for any $p>0$.

Finally consider the first integral in \eqref{intermediate1}
\begin{eqnarray*}
&&\int_{-\delta\sqrt{n}}^{\delta\sqrt{n}}|e^{n\psi(\theta/\sqrt{n})}-e^{-\theta^2/2}q_n(\theta)|\frac{d\theta}{2\pi}\\
&\leq&\int_{-\delta\sqrt{n}}^{\delta\sqrt{n}}e^{-\theta^2/2}|e^{n\psi(\theta/\sqrt{n})+\theta^2/2}-w_n(\theta)|\frac{d\theta}{2\pi}+\int_{-\delta\sqrt{n}}^{\delta\sqrt{n}}e^{-\theta^2/2}|r_n(\theta)|\frac{d\theta}{2\pi}
\end{eqnarray*}
By our observation on $r_n(\theta)$ that it can be written as a sum of $P_k(\theta,n)n^{-k}$ with conditions on $P_k$ and $k$, together with the fact that $\int_{-\delta\sqrt{n}}^{\delta\sqrt{n}}e^{-\theta^2/2}|\theta|^bd\theta<\infty$ for any $b>0$ and the slowly varying functions $L$ and $\zeta_L$ grow slower than any polynomial, we can verify that
$$\int_{-\delta\sqrt{n}}^{\delta\sqrt{n}}e^{-\theta^2/2}|r_n(\theta)|\frac{d\theta}{2\pi}=\left\{\begin{array}{ll}o\left(\frac{L(n^{1/2})}{n^{\alpha/2-1}}\right)&\mbox{\ for non-even $\alpha$}\\
o\left(\frac{\zeta_L(n^{1/2})}{n^{\alpha/2-1}}\right)&\mbox{\ for even
$\alpha$}\end{array}\right.$$ Now consider the first part. We use
the estimate
\begin{equation}
\left|e^\alpha-\sum_{j=0}^\gamma\frac{\beta^j}{j!}\right|\leq|e^\alpha-e^\beta|+\left|e^\beta-\sum_{j=0}^\gamma\frac{\beta^j}{j!}\right|\leq\left(|\alpha-\beta|+\frac{|\beta|^{\gamma+1}}{(\gamma+1)!}\right)e^{\max\{\alpha,\beta\}} \label{intermediate2}
\end{equation}
where we use mean value theorem in the last inequality. Also, we can
find $\delta$ small enough so that for $\theta\in[-\delta,\delta]$,
we have
$$\left|\psi(\theta)+\frac{\theta^2}{2}\right|\leq\frac{\theta^2}{4} \mbox{\ \ and\ \ }|\lambda(\theta)|\leq\frac{\theta^2}{4}$$
which implies that for $\theta\in[-\delta\sqrt{n},\delta\sqrt{n}]$,
we get
\begin{equation}
\left|n\psi\left(\frac{\theta}{\sqrt{n}}\right)+\frac{\theta^2}{2}\right|\leq\frac{\theta^2}{4} \mbox{\ \ and\ \ }\left|n\lambda\left(\frac{\theta}{\sqrt{n}}\right)\right|\leq\frac{\theta^2}{4} \label{intermediate3}
\end{equation}
Using the above estimates \eqref{intermediate2} and \eqref{intermediate3},
\begin{eqnarray}
&&\int_{-\delta\sqrt{n}}^{\delta\sqrt{n}}e^{-\theta^2/2}|e^{n\psi(\theta/\sqrt{n})+\theta^2/2}-w_n(\theta)|\frac{d\theta}{2\pi} \notag\\
&\leq&\int_{-\delta\sqrt{n}}^{\delta\sqrt{n}}e^{-\theta^2/4}\left|nR\left(\frac{\theta}{\sqrt{n}}\right)\right|\frac{d\theta}{2\pi}+\int_{-\delta\sqrt{n}}^{\delta\sqrt{n}}e^{-\theta^2/4}\left|\frac{n\lambda(\theta/\sqrt{n})}{(\lfloor\alpha/2-1\rfloor+1)!}\right|^{\lfloor\alpha/2-1\rfloor+1}\frac{d\theta}{2\pi} \label{intermediate4}
\end{eqnarray}
Now since
$$\lambda(\theta)=o(\theta^{2+p})$$
for $\theta\in[-\delta\sqrt{n},\delta\sqrt{n}]$ and some small $p>0$, it is easy to see that the second integral in \eqref{intermediate4}
$$\int_{-\delta\sqrt{n}}^{\delta\sqrt{n}}e^{-\theta^2/4}\left|\frac{n\lambda(\theta/\sqrt{n})}{(\lfloor\alpha/2-1\rfloor+1)!}\right|^{\lfloor\alpha/2-1\rfloor+1}\frac{d\theta}{2\pi}=o\left(\frac{1}{n^{\alpha/2-1}}\right)$$
On the other hand, we have
\begin{equation}
R(\theta)=\left\{\begin{array}{ll}o\left(|\theta|^\alpha
L\left(\frac{1}{|\theta|}\right)\right)&\mbox{\ for non-even
$\alpha$}\\
o\left(|\theta|^\alpha\zeta_L\left(\frac{1}{|\theta|}\right)\right)&\mbox{\
for even $\alpha$}\end{array}\right. \label{R}
\end{equation} Suppose that $\alpha$ is
non-even. For given $\epsilon$, we can choose $\delta$ small enough
such that for $\theta\in[-\delta,\delta]$,
$$R(\theta)\leq\epsilon|\theta|^\alpha
L\left(\frac{1}{|\theta|}\right)$$ which implies that for
$\theta\in[-\delta\sqrt{n},\delta\sqrt{n}]$,
$$nR\left(\frac{\theta}{\sqrt{n}}\right)\leq\epsilon\frac{|\theta|^\alpha}{n^{\alpha/2-1}}L\left(\frac{n^{1/2}}{|\theta|}\right)$$
and the first integral in \eqref{intermediate4}
$$\int_{-\delta\sqrt{n}}^{\delta\sqrt{n}}e^{-\theta^2/4}\left|nR\left(\frac{\theta}{\sqrt{n}}\right)\right|\frac{d\theta}{2\pi}\leq\epsilon\int_{-\delta\sqrt{n}}^{\delta\sqrt{n}}e^{-\theta^2/4}\frac{|\theta|^\alpha}{n^{\alpha/2-1}}L\left(\frac{n^{1/2}}{|\theta|}\right)\frac{d\theta}{2\pi}$$
Note that by Lemma \ref{regularlyvaryinglemma} we can pick $M>0$ and
any $\rho_1,\rho_2>0$ such that
\begin{equation*}
\frac{L(\frac{n^{1/2}}{|\theta|})}{L(n^{1/2})}\leq\left\{%
\begin{array}{ll}
C|\theta|^{-\rho_1} & \text{for\ }|\theta|<1 \\
C|\theta|^{\rho_2} & \text{for\ }|\theta|\geq 1%
\end{array}%
\right.
\end{equation*}
for $n>M$. So by dominated convergence theorem and that
\begin{equation}
\frac{L(\frac{n^{1/2}}{|\theta|})}{L(n^{1/2})}\to1
\label{convergence}
\end{equation}
as $n\to\infty$ for all fixed nonzero $\theta$, we have
\begin{eqnarray*}
&&\int_{-\delta\sqrt{n}}^{\delta\sqrt{n}}e^{-\theta^2/4}\frac{|\theta|^\alpha}{n^{\alpha/2-1}}L\left(\frac{n^{1/2}}{|\theta|}\right)\frac{d\theta}{2\pi}\\
&=&\int_{-\delta\sqrt{n}}^{\delta\sqrt{n}}e^{-\theta^2/4}|\theta|^\alpha\frac{L(n^{1/2}/|\theta|)}{L(n^{1/2})}\frac{d\theta}{2\pi}\frac{L(n^{1/2})}{n^{\alpha/2-1}}\\
&\leq&\int_{-\infty}^\infty e^{-\theta^2/4}|\theta|^\alpha\frac{L(n^{1/2}/|\theta|)}{L(n^{1/2})}\frac{d\theta}{2\pi}\frac{L(n^{1/2})}{n^{\alpha/2-1}}\\
&\sim&\int_{-\infty}^\infty
e^{-\theta^2/4}|\theta|^\alpha\frac{d\theta}{2\pi}\frac{L(n^{1/2})}{n^{\alpha/2-1}}
\end{eqnarray*}
as $n\nearrow\infty$. Hence
$$\limsup_{n\to\infty}\frac{\int_{-\delta\sqrt{n}}^{\delta\sqrt{n}}e^{-\theta^2/4}\left|nR\left(\frac{\theta}{\sqrt{n}}\right)\right|\frac{d\theta}{2\pi}}{L(n^{1/2})/n^{\alpha/2-1}}\leq\epsilon C$$
for some constant $C>0$. Since $\epsilon$ is arbitrary, we get
$$\int_{-\delta\sqrt{n}}^{\delta\sqrt{n}}e^{-\theta^2/4}\left|nR\left(\frac{\theta}{\sqrt{n}}\right)\right|\frac{d\theta}{2\pi}=o\left(\frac{L(n^{1/2})}{n^{\alpha/2-1}}\right)=o\left(\frac{L(n^{1/2}x)}{n^{\alpha/2-1}}\right)$$
for any fixed  $x$. The case for even $\alpha$ is the same except
$L$ is replaced by $\zeta_L$ that is also a slowly varying function.

We can therefore conclude that
$$\int_{-\infty}^\infty\left[e^{-i\theta
x+n\psi(\theta/\sqrt{n})}-e^{-i\theta
x-\theta^2/2}q_n(\theta)\right]\frac{d\theta}{2\pi}=\left\{\begin{array}{ll}o\left(\frac{L(n^{1/2}x)}{n^{\alpha/2-1}}\right)&\mbox{\
for non-even
$\alpha$}\\
o\left(\frac{\zeta_L(n^{1/2}x)}{n^{\alpha/2-1}}\right)&\mbox{\ for
even $\alpha$}\end{array}\right.$$

We now evaluate $\int_{-\infty}^\infty e^{-i\theta
x-\theta^2/2}q_n(\theta)\frac{d\theta}{2\pi}$. Note that
$$\int_{-\infty}^\infty e^{-i\theta
x-\theta^2/2}q_n(\theta)\frac{d\theta}{2\pi}=\eta(x)\left[1+\sum_{\substack{ 3\leq j<\alpha  \\ j\text{\ even}}}%
\frac{G_j(x)}{n^{j/2-1}}\right]+\int_{-\infty}^\infty e^{-i\theta
x-\theta^2/2}n\xi\left(\frac{\theta}{n^{1/2}}\right)\frac{d\theta}{2\pi}$$
where $G_j$'s are the ordinary Edgeworth coefficients. So it suffices to focus on the last term
\begin{eqnarray*}
&&\int_{-\infty}^\infty e^{-i\theta
x-\theta^2/2}n\xi\left(\frac{\theta}{n^{1/2}}\right)\frac{d\theta}{2\pi}\\
&=&\left\{\begin{array}{ll} -\frac{\pi}{\Gamma(\alpha+1)\sin
\frac{\alpha\pi}{2}}\int_{-\infty}^\infty e^{-i\theta x-\theta^2/2}
|\theta|^\alpha
L(\frac{n^{1/2}}{|\theta|})\frac{d\theta}{2\pi}\ \frac{1}{n^{\alpha/2-1}}&\mbox{\ for non-even $\alpha$}\\
\frac{2(-1)^{\alpha/2}}{\Gamma(\alpha+1)}\ \int_{-\infty}^\infty
e^{-i\theta x-\theta^2/2}|\theta|^\alpha
\zeta_L(\frac{n^{1/2}}{|\theta|})\frac{d\theta}{2\pi}\
\frac{1}{n^{\alpha/2-1}}&\mbox{\ for even
$\alpha$}\end{array}\right.\end{eqnarray*} By dominated convergence
theorem and (\ref{convergence}), we have for odd $\alpha$,
\begin{eqnarray*}
&&-\frac{\pi}{\Gamma(\alpha+1)\sin
\frac{\alpha\pi}{2}}\int_{-\infty}^\infty e^{-i\theta x-\theta^2/2}
|\theta|^\alpha
L\left(\frac{n^{1/2}}{|\theta|}\right)\frac{d\theta}{2\pi}\
\frac{1}{n^{\alpha/2-1}}\\
&=&-\frac{\pi}{\Gamma(\alpha+1)\sin
\frac{\alpha\pi}{2}}\int_{-\infty}^\infty e^{-i\theta x-\theta^2/2}
|\theta|^\alpha\frac{L(\frac{n^{1/2}}{|\theta|})}{L(n^{1/2})}\frac{d\theta}{2\pi}\
\frac{L(n^{1/2})}{n^{\alpha/2-1}}\\
&\sim&-\frac{\pi}{\Gamma(\alpha+1)\sin
\frac{\alpha\pi}{2}}\int_{-\infty}^\infty e^{-i\theta x-\theta^2/2}
|\theta|^\alpha \frac{d\theta}{2\pi}\
\frac{L(n^{1/2}x)}{n^{\alpha/2-1}}\\
&=&-\frac{\sqrt{2}}{\alpha!}\frac{d^\alpha}{dx^\alpha}D(\frac{x}{\sqrt{2}})\frac{L(n^{1/2}x)}{n^{\alpha/2-1}}
\end{eqnarray*}
by using Lemma \ref{Edgeworthlemma} in the last equality.

For non-integral and even $\alpha$ essentially the same argument
follows, again noting that $\zeta_L(x)$ is also a slowly varying
function for proving the even $\alpha$ case.
\end{proof}
\bigskip

\section{Beyond Central Limit Region}

Although this paper concerns expansions in the central limit region,
we devote this section to an interesting connection of our extra
term with the large deviations of regularly varying rv's.

By local analysis (see for example, Bender and Orszag (1999))
\begin{equation*}
\frac{d^\alpha}{dx^\alpha}D(\frac{x}{\sqrt{2}})\sim(-1)^\alpha\frac{\alpha!}{%
\sqrt{2}x^{1+\alpha}}
\end{equation*}
as $x\to\infty$ and for non-integral $\nu$%
,
\begin{equation*}
D_\nu(z)\sim\left\{%
\begin{array}{ll}
z^\nu e^{-z^2/4} & \text{as\ }z\to+\infty \\
\frac{\sqrt{2\pi}}{\Gamma(-\nu)}e^{z^2/4}|z|^{-\nu-1} & \text{as\ }%
z\to-\infty%
\end{array}%
\right.
\end{equation*}

So for odd and non-integral $\alpha$,
\begin{equation*}
G_\alpha(x)\sim x^{-(1+\alpha)}
\end{equation*}
as $x\to\infty$. Note that the result for non-integral $\alpha$ follows from
the observation that the negative side of $D_\alpha$ in
$G_\alpha(x)$ dominates and the reflection property of gamma
function. Thus when $\alpha$ is non-even, we get
\begin{equation}
G_\alpha(x)\frac{L(n^{1/2}x)}{n^{\alpha/2-1}}\sim
\frac{L(n^{1/2}x)}{x^{1+\alpha}n^{\alpha/2-1}} \label{largedeviation}
\end{equation}
for $x\to\infty$.

On the other hand, Rozovskii (1989, 1993) proved a result for regularly
varying rv on the whole axis:
\begin{equation}
P(S_n>x)=\left[\bar{\Phi}\left(\frac{x}{\sqrt{n}}\right)+nP(X_k>x)\right](1+o(1))
\label{realaxisasymptotic} \end{equation} uniformly for
$x>\sqrt{n}$, where $\bar{\Phi}(\cdot)$ denotes the complementary normal distribution. In other words,
$$\lim_{n\nearrow\infty}\sup_{x>\sqrt{n}}\left|\frac{P(S_n>x)}{\bar{\Phi}(x/\sqrt{n})+nP(X_k>x)}-1\right|=0$$

In particular, the result implies that for $x>\sqrt{(\alpha-2)n\log
n}$, the large deviations asymptotic (the second term on the right
hand side) dominates i.e. $P(S_n>x)\sim nP(X_k>x)$. If we assume the
existence of density, this can be written as
$$f_{S_n/\sqrt{n}}(x)\sim nf(x\sqrt{n})\cdot\sqrt{n}$$
for $x>\sqrt{(\alpha-2)\log n}$, which is exactly
(\ref{largedeviation}). In other words, the non-analytic term
matches exactly with the large deviations asymptotic.

However, for even $\alpha$ our extra term resembles standard
Edgeworth expansion, and there is no region on the real axis that
this can coincide with the large deviation asymptotic.

To rigorously show the validity of the connection would require
knowledge about the second-order variation of the slowly varying
function $L$ i.e. how fast $L(xt)/L(x)\to1$ for given $t$. This will
be a direction for future research, but in this paper we are
satisfied that at least in special cases this connection is easily
seen to be valid:

\begin{thm}
Assume $(X_k:k\geq1)$ are iid rv's with unit variance and symmetric
Pareto density $f(x)=a_f/(b_f+|x|^{\alpha+1})$ where $\alpha$ is
non-even, and $a_f$ and $b_f$ are constants. Then we have the
asymptotic
\begin{equation}
f_{S_n/n^{1/2}}(x)=\eta(x)\left[1+\sum_{\substack{ 3\leq j< \alpha  \\ j\text{\ even}}}%
\frac{G_j(x)}{n^{j/2-1}}\right]+\frac{G_\alpha(x)}{n^{\alpha/2-1}}+o\left(\frac{G_\alpha(x)}{n^{\alpha/2-1}}\right)
\label{refinedrealaxisasymptotic}
\end{equation}
uniformly over $x\leq\sqrt{(\alpha-2)\log n}$, where $G_j(x)$ are defined as in
Theorem \ref{symmetrictheorem}.
\label{moderatedeviations}
\end{thm}

In other words, the refined moderate deviations over the range
$x\leq\sqrt{(\alpha-2)\log n}$ is the same as the Edgeworth expansion together with our additional term.

\begin{proof}
The proof is a simple adaptation from the proof of Theorem
\ref{symmetrictheorem}. Note that in that proof we bound uniformly
the error of Edgeworth expansion over the real axis, but this
error can be too coarse for $x$ outside the central limit region. It
suffices, therefore, to refine our error estimate so that it always
has smaller order than our Edgeworth expansion over
$x\leq\sqrt{(\alpha-2)\log n}$. In other words, we need to show that the error of our approximation using Edgeworth expansion is of order $o(L(n^{1/2}x)/(x^{1+\alpha}n^{\alpha/2-1}))$ for any $x\leq\sqrt{(\alpha-2)\log n}$.

We will first prove the result for the odd $\alpha$ case. An easy examination of our proof of
Theorem \ref{symmetrictheorem} would reveal that for Pareto rv, the
only error that needs further refinement is that caused by
$$\int_{-\delta\sqrt{n}}^{\delta\sqrt{n}}e^{-\theta^2/4}\left|nR\left(\frac{\theta}{\sqrt{n}}\right)\right|\frac{d\theta}{2\pi}$$
For this we need a finer estimate of $R(\theta)$, which in turn
needs a finer version of Case 1 in Lemma \ref{paretolemma}. Note
that
\begin{align*}
\int_{-\infty}^\infty\frac{a_fx^{\alpha-1}(e^{i\theta
x}-1)}{b_f+|x|^{1+\alpha}}dx&=2a_f\int_0^\infty\frac{x^{\alpha-1}(\cos(\theta
x)-1)}{b_f+x^{1+\alpha}}dx\\
&=2a_f\theta\int_0^\infty\frac{u^{\alpha-1}(\cos
u-1)}{b_f\theta^{1+\alpha}+u^{1+\alpha}}du\\
&=2a_f\theta\int_0^\infty\frac{\cos
u-1}{u^2}du-2a_fb_f\theta^{2+\alpha}\int_0^\infty\frac{(\cos
u-1)}{u^2(b_f\theta^{1+\alpha}+u^{1+\alpha})}du
\end{align*}
The first part becomes $-a_f\pi\theta$, while the second part is
\begin{eqnarray*}
&&-2a_fb_f\theta^{2+\alpha}\int_0^\infty\frac{(\cos
u-1)}{u^2(b_f\theta^{1+\alpha}+u^{1+\alpha})}du\\
&=&-2a_fb_f\theta^{2+\alpha}\left[\int_0^\epsilon\frac{(\cos
u-1)}{u^2(b_f\theta^{1+\alpha}+u^{1+\alpha})}du+\int_\epsilon^\infty\frac{(\cos
u-1)}{u^2(b_f\theta^{1+\alpha}+u^{1+\alpha})}du\right]
\end{eqnarray*}
where $\epsilon$ is a small number. Note that for $0<u<\epsilon$ we
have $(\cos u-1)/u^2$ bounded around $-1$. Hence
\begin{align*}
\left|\int_0^\epsilon\frac{(\cos
u-1)}{u^2(b_f\theta^{1+\alpha}+u^{1+\alpha})}du\right|&\leq
C\int_0^\epsilon\frac{1}{b_f\theta^{1+\alpha}+u^{1+\alpha}}du\\
&=\frac{C}{\theta^\alpha}\int_0^{\frac{\epsilon}{\theta}}\frac{1}{b_f+x^{1+\alpha}}dx\\
&\leq\frac{C}{\theta^\alpha}
\end{align*}
where $C$ are positive constants. Obviously
$$\int_\epsilon^\infty\frac{(\cos
u-1)}{u^2(b_f\theta^{1+\alpha}+u^{1+\alpha})}du=O(1)$$ Therefore
$$\int_{-\infty}^\infty\frac{a_fx^{\alpha-1}(e^{i\theta
x}-1)}{b_f+|x|^{1+\alpha}}dx=-2a_f\pi\theta+O(\theta^2)$$ and hence
$R(\theta)=O(|\theta|^{1+\alpha})$. Replacing \eqref{R} with this estimate of $R(\theta)$ and applying the same argument as in the proof of Theorem 2 implies that
$$\int_{-\delta\sqrt{n}}^{\delta\sqrt{n}}e^{-\theta^2/4}\left|nR\left(\frac{\theta}{\sqrt{n}}\right)\right|\frac{d\theta}{2\pi}=O\left(\frac{1}{n^{\alpha/2-1/2}}\right)=o\left(\frac{L(n^{1/2}x)}{x^{1+\alpha}n^{\alpha/2-1}}\right)$$
for $x\leq\sqrt{(\alpha-2)\log n}$. This leads to a uniform estimate of our error in
the Edgeworth expansion that always has order less than
$G_\alpha(x)/n^{\alpha/2-1}$ over $x\leq\sqrt{(\alpha-2)\log n}$.

Now for non-integer $\alpha$, when $q=\lfloor\alpha\rfloor$ is even, Case 3 in Lemma \ref{paretolemma} can be written as
\begin{align*}
\int_{-\infty}^\infty\frac{a_fx^q(e^{i\theta x}-1)}{b_f+|x|^{1+\alpha}}dx&=2a_f\theta^{\alpha-q}\int_0^\infty\frac{u^q}{b_f\theta^{1+\alpha}+u^{1+\alpha}}(\cos u-1)du\\
&=2a_f\left[\theta^{\alpha-q}\int_0^\infty\frac{\cos u-1}{u^{1+\alpha-q}}du-b_f\theta^{1+\alpha}\int_0^\infty\frac{\cos u-1}{u^{1+\alpha-q}(b_f\theta^{1+\alpha}+u^{1+\alpha})}du\right]
\end{align*}
Now 
$$\int_0^\infty\frac{\cos u-1}{u^{1+\alpha-q}(b_f\theta^{1+\alpha}+u^{1+\alpha})}du=\int_0^\epsilon\frac{\cos u-1}{u^{1+\alpha-q}(b_f\theta^{1+\alpha}+u^{1+\alpha})}du+\int_\epsilon^\infty\frac{\cos u-1}{u^{1+\alpha-q}(b_f\theta^{1+\alpha}+u^{1+\alpha})}du$$
for small $\epsilon>0$. Clearly the second integral is bounded. Consider the first integral
\begin{align*}
\left|\int_0^\epsilon\frac{\cos u-1}{u^{1+\alpha-q}(b_f\theta^{1+\alpha}+u^{1+\alpha})}du\right|&\leq C\int_0^\epsilon\frac{u^{1-(\alpha-q)}}{b_f\theta^{1+\alpha}+u^{1+\alpha}}du\\
&=\frac{C}{\theta^{\alpha-1+(\alpha-q)}}\int_0^{\frac{\epsilon}{\theta}}\frac{y^{1-(\alpha-q)}}{b_f+y^{1+\alpha}}dy\\
&=O\left(\frac{1}{\theta^{\alpha-1+(\alpha-q)}}\right)
\end{align*}
where $C$ is a positive constant, and a substitution $u=\theta y$ is used in the first equality. Hence
$$\int_{-\infty}^\infty\frac{a_fx^q(e^{i\theta x}-1)}{b_f+|x|^{1+\alpha}}dx=2a_f\theta^{\alpha-q}\Gamma(-\alpha+q)\cos\left(\frac{-\alpha+q}{2}\pi\right)+O(\theta^2)$$

When $q=\lfloor\alpha\rfloor$ is odd, we have
\begin{align*}
\int_{-\infty}^\infty\frac{a_fx^qe^{i\theta x}}{b_f+|x|^{1+\alpha}}dx&=2ia_f\theta^{\alpha-q}\int_0^\infty\frac{u^q\sin u}{b_f\theta^{1+\alpha}+u^{1+\alpha}}du\\
&=2ia_f\left[\theta^{\alpha-q}\int_0^\infty\frac{\sin u}{u^{1+\alpha-q}}du-b_f\theta^{1+\alpha}\int_0^\infty\frac{\sin u}{u^{1+\alpha-q}(b_f\theta^{1+\alpha}+u^{1+\alpha})}du\right]
\end{align*}
Now
$$\int_0^\infty\frac{\sin u}{u^{1+\alpha-q}(b_f\theta^{1+\alpha}+u^{1+\alpha})}du=\int_0^\epsilon\frac{\sin u}{u^{1+\alpha-q}(b_f\theta^{1+\alpha}+u^{1+\alpha})}du+\int_\epsilon^\infty\frac{\sin u}{u^{1+\alpha-q}(b_f\theta^{1+\alpha}+u^{1+\alpha})}du$$
for small $\epsilon>0$. Clearly the second integral is bounded. Consider the first integral
\begin{align*}
\left|\int_0^\epsilon\frac{\sin u}{u^{1+\alpha-q}(b_f\theta^{1+\alpha}+u^{1+\alpha})}du\right|&\leq C\int_0^\epsilon\frac{1}{u^{\alpha-q}(b_f\theta^{1+\alpha}+u^{1+\alpha})}du\\
&=\frac{C}{\theta^{\alpha+(\alpha-q)}}\int_0^{\frac{\epsilon}{\theta}}\frac{1}{y^{\alpha-q}(b_f+y^{1+\alpha})}dy\\
&=O\left(\frac{1}{\theta^{\alpha+(\alpha-q)}}\right)
\end{align*}
where $C$ is a positive constant, and a substitution $u=\theta y$ is used in the first equality. Hence
$$\int_{-\infty}^\infty\frac{a_fx^qe^{i\theta x}}{b_f+|x|^{1+\alpha}}dx=2ia_f\theta^{\alpha-q}\Gamma(-\alpha+q)\sin\left(\frac{-\alpha+q}{2}\pi\right)+O(\theta)$$

Therefore
$$R(\theta)=\left\{\begin{array}{ll}O(|\theta|^{\ulcorner\alpha\urcorner+1})&\mbox{\
when $\llcorner\alpha\lrcorner$ is even}\\
O(|\theta|^{\ulcorner\alpha\urcorner})&\mbox{\ when
$\llcorner\alpha\lrcorner$ is odd}\end{array}\right.$$ and the same
argument as the odd $\alpha$ case holds.

\end{proof}

\section{Non-Symmetric Density}

With the result for symmetric density in hand, we can generalize our
results to the non-symmetric counterpart, namely Theorem
\ref{maintheorem}. The main idea is to write a non-symmetric density
as
\begin{equation*}
f(x)=r(x)+s(x)
\end{equation*}
where
\begin{align}
r(x)=\frac{1}{2}[f(x)+f(-x)]\label{evenpart}
\end{align}
is a symmetric function, and
\begin{align}
s(x)=\frac{1}{2}[f(x)-f(-x)]\label{oddpart}
\end{align}
is an odd function. The symmetric function $r$ can be easily handled
using the results in the last section, and the odd function $s$ can
also be dealt with by slight modification. Indeed we have the
following proposition:

\begin{proposition}
For rv following (\ref{density}) with $L_+$ and $L_-$ satisfying
Assumption \ref{mainassumption}, we can decompose its density as
\begin{equation*}
f(x)=r(x)+s(x)
\end{equation*}
where $r$ and $s$ are defined by (\ref{evenpart}) and
(\ref{oddpart}) respectively. Then $r$ and $s$ must be regularly
varying and have the same order $\alpha\equiv\min(\beta,\gamma)$,
and we can write
\begin{equation}
r(x)=\frac{L_r(x)}{1+|x|^{1+\alpha}} \ \
s(x)=\frac{L_s(x)}{1+|x|^{1+\alpha}} \ \ ,\ -\infty<x<\infty
\label{rsdefinition}
\end{equation}
where $L_r$ is even and $L_s$ is odd, and they are both slowly
varying functions satisfying Assumption \ref{mainassumption}. Also
define
\begin{equation*}
m_n^{(r)}=\int_{-\infty}^\infty x^nr(x)dx \ \
m_n^{(s)}=\int_{-\infty}^\infty x^ns(x)dx \text{ \ \ \ for \ $0\leq
n\leq\alpha-1$}
\end{equation*}
as the ``moments" of $r$ and $s$. The cumulant generating function
takes the form
$$\psi(\theta)=\chi(\theta)+\xi(\theta)+o(\xi(\theta))$$
where
$$\chi(\theta)=1+\sum_{1\leq j<\alpha}\kappa_j\frac{(i\theta)^j}{j!}$$ is the ordinary Taylor series
expansion up to the largest cumulant, and
\begin{equation}
\xi(\theta)=\left\{\begin{array}{ll}\frac{2i^\alpha}{\Gamma(\alpha+1)}\
|\theta|^\alpha\left[\zeta_{L_r}(\frac{1}{|\theta|})\pm i\frac{\pi}{2}L_s(\frac{1}{|\theta|})\right]&\mbox{\
for even
$\alpha$}\\
\frac{2i^\alpha}{\Gamma(\alpha+1)}\
|\theta|^\alpha\left[\pm\zeta_{L_s}(\frac{1}{|\theta|})+i\frac{\pi}{2}L_r(\frac{1}{|\theta|})\right]&\mbox{\
for odd
$\alpha$}\\
-\frac{2\pi}{\Gamma(\alpha+1)\sin(\alpha\pi)}\
|\theta|^\alpha\left[\cos\frac{\alpha\pi}{2}L_r(\frac{1}{|\theta|})\mp i\sin\frac{\alpha\pi}{2}L_s(\frac{1}{|\theta|})\right]&\mbox{\
for non-integer $\alpha$} \end{array}\right. \label{nonsymmetricxi}
\end{equation} is the non-analytic component. Here $\pm$ and $\mp$ both refer to the cases when $\theta>0$ or $\theta<0$.
\label{nonsymmetriccumulant}
\end{proposition}

\bigskip

\begin{proof} That $r$ and $s$ are regularly varying, have the same order and
can be written as (\ref{rsdefinition}) satisfying Assumption
\ref{mainassumption} is obvious. Now define $\phi_r(\theta)$ and
$\phi_s(\theta)$ as the Fourier transforms, or ``characteristic
functions" of $r$ and $s$. Since $r$ is symmetric, by the same
argument as Proposition \ref{regularlyvaryingcharacteristic}, we get
\begin{equation}
\phi_r(\theta) = 1+\sum_{\substack{ 2\leq j<\alpha  \\ j \text{\
even}}}m_j^{(r)}\frac{(i\theta)^j}{j!}+
\left\{\begin{array}{ll}-\frac{\pi }{\Gamma(\alpha+1)\sin
\frac{\alpha\pi}{2}}\ |\theta|^\alpha
L_r\left(\frac{1}{|\theta|}\right)+o\left(|\theta|^\alpha
L_r\left(\frac{1}{|\theta|}\right)\right)&\mbox{\ for non-even
$\alpha$}\\\frac{2(-1)^{\alpha/2}}{\Gamma(\alpha+1)}\
|\theta|^\alpha\zeta_{L_r}\left(\frac{1}{|\theta|}\right)+o\left(|\theta|^\alpha\zeta_{L_r}\left(\frac{1}{|\theta|}\right)\right)&\mbox{\ for
even $\alpha$}\end{array}\right. \label{rcharacteristic}
\end{equation}
Now consider the odd function $s$. Note that the role played by odd
and even $\alpha$ (or $q$) is reversed
because of odd instead of even function. Also note that $\phi_s%
(\theta)=-\phi_s(-\theta)$, so it suffices to prove the result for
$\theta>0$ and the $\theta<0$ case follows. Following similar line
of proof as Proposition \ref{regularlyvaryingcharacteristic}, we get
\begin{equation}
\phi_s(\theta) = 1+\sum_{\substack{ 1\leq j<\alpha  \\ j \text{\
odd}}}m_j^{(s)}\frac{(i\theta)^j}{j!}+
\left\{\begin{array}{ll}\frac{i\pi }{\Gamma(\alpha+1)\cos
\frac{\alpha\pi}{2}}\ |\theta|^\alpha
L_s\left(\frac{1}{\theta}\right)+o\left(|\theta|^\alpha
L_s\left(\frac{1}{\theta}\right)\right)&\mbox{\ for non-odd
$\alpha$}\\\frac{2i^\alpha}{\Gamma(\alpha+1)}\
|\theta|^\alpha\zeta_{L_s}\left(\frac{1}{\theta}\right)+o\left(|\theta|^\alpha\zeta_{L_s}\left(\frac{1}{\theta}\right)\right)&\mbox{\ for
odd $\alpha$}\end{array}\right. \label{scharacteristic}
\end{equation}
for $\theta>0$. Adding (\ref{rcharacteristic}) and (\ref{scharacteristic}), and
noting that the representation can be carried over to cumulant
generating function, we obtain the stated result.

\end{proof}

\bigskip

We also need a correspondence of Lemma \ref{Edgeworthlemma} for odd
function, the proof of which is very similar and is thus omitted:

\begin{lemma}
\begin{equation}
\int_{-\infty}^\infty e^{-i\theta x-\theta^2/2}(\pm|\theta|^\alpha)
d\theta=\left\{\begin{array}{ll}2\sqrt{2}i^{\alpha+1}\frac{d^\alpha}{dx^\alpha}D(\frac{x}{\sqrt{2}})&\mbox{\
for even $\alpha$}\\
-\sqrt{2\pi}e^{-x^2/2}i^\alpha H_\alpha(x)&\mbox{\ for odd
$\alpha$}\\
i\sqrt{\frac{\pi}{2}}\csc\frac{\alpha\pi}{2}e^{-x^2/4}[D_\alpha(x)-D_\alpha(-x)]&\mbox{\
for non-integer $\alpha$}
\end{array}\right.
\end{equation}
where $D(z)$, $H_k(z)$ and $D_\nu(z)$ are Dawson integral, Hermite
polynomial of order $k$ and classical parabolic cylinder function
with parameter $\nu$ respectively. Here
$$\pm|\theta|^\alpha=\left\{\begin{array}{ll}|\theta|^\alpha&\mbox{for $\theta>0$}\\
-|\theta|^\alpha&\mbox{for $\theta<0$}
\end{array}\right.$$

\label{Edgeworthoddlemma}\end{lemma}

\bigskip

The proof of Theorem \ref{maintheorem} follows easily from
Proposition \ref{nonsymmetriccumulant} and Lemma
\ref{Edgeworthoddlemma}:

\begin{proof}[Proof of Theorem \ref{maintheorem}] By (\ref{density}) we can write $L_r$ and $L_s$ as
\begin{equation}
\begin{array}{lll} L_r(x)=\frac{1}{2}[L_+(|x|)+L_-(|x|)]&
L_s(x)=\pm\frac{1}{2}[L_+(|x|)-L_-(|x|)]&\mbox{\ when
$\beta=\gamma$}\\
L_r(x)\sim\frac{1}{2}L_+(|x|)&
L_s(x)\sim\pm\frac{1}{2}L_+(|x|)&\mbox{\ when $\beta<\gamma$}\\
L_r(x)\sim\frac{1}{2}L_-(|x|)&
L_s(x)\sim\mp\frac{1}{2}L_-(|x|)&\mbox{\ when $\beta>\gamma$}
\end{array} \label{nonsymmetricfunctions}
\end{equation}
where the $\pm$ and $\mp$ signs correspond to $x\geq0$ and $x<0$. We
can proceed similarly as in the proof of Theorem
\ref{symmetrictheorem} by substituting (\ref{nonsymmetricfunctions})
into $L_r$ and $L_s$, and using Lemma \ref{Edgeworthoddlemma}. Note
the use of (\ref{integrability}) in obtaining the $\beta<\gamma$ and
$\beta>\gamma$ cases.
\end{proof}

\bigskip

\section{References}

\begin{enumerate}
\item Bazant, M. Z. (2006), \emph{18.366 Random Walks and Diffusion}, Lecture Notes, MIT OpenCourseWare,
http://ocw.mit.edu/OcwWeb/Mathematics/18-366Fall-2006/LectureNotes/index.htm (See the notes and supplementary material for lecture 7.)

\item Bender, C. M., and Orszag, S. A. (1999), \emph{Advanced Mathematical
Methods for Scientists and Engineers: Asymptotic Methods and
Perturbation Theory}, 2nd edition, Springer-Verlag.

\item Bouchaud, J. P. and Georges, M. (1990), Anomalous diffusion in disordered media: Statistical mechanisms, models, and physical applications, \emph{Physics Reports} {\bf 195} (4-5), 127-293.

\item Bouchaud, J. P., and Potters, M. (2000), \emph{Theory of Financial Risks}, Cambridge, UK: Cambridge University Press.

\item Embrechts, P., Kluppelberg C., and Mikosch, T. (1997),
\emph{Modelling Extremal Events for Insurance and Finance},
Springer.

\item Evans, M., and Swartz, T. (2000), \emph{Approximating
Integrals via Monte Carlo and Deterministic Methods}, Oxford
Statistical Science Series.

\item Field, C. A., and Ronchetti, E. (1990), \emph{Small Sample
Asymptotics}, Institute of Mathematical Statistics, Lecture Notes -
Monograph Series.

\item Hall, P. (1992), \emph{The Bootstrap and Edgeworth Expansion},
Springer-Verlag.

\item Hughes, B. (1996), \emph{Random Walks and Random Environments}, Vol. 1, Oxford, UK: Clarendon Press.

\item Metzler, R. and Klafter, J. (2000), The random walk's guide to anomalous diffusion: a fractional dynamics approach, \emph{Physics Reports} {\bf 339}, 1-77.

\item Resnick, S. I. (1987), \emph{Extreme Values, Regular Variation,
and Point Processes}, Springer-Verlag.

\item Rozovskii, L. V. (1989), Probabilities of large deviations of
sums of independent random variables with common distribution
function in the domain of attraction of the normal law, \emph{Theory
Probab. Appl.} \textbf{34}(4), 625-644.

\item Vinogradov, V. (1994), \emph{Refined Large Deviation Limit
Theorems}, Longman Scientific \& Technical.

\item Weiss, G. H. (1994), \emph{ Aspects and Applications of the Random Walk }, North-Holland.

\end{enumerate}

\end{document}